%% file: main-arxiv.tex
\documentclass{ar-1col-S2O}
\usepackage[numbers]{natbib}
\usepackage{url}
\usepackage{amsmath,amssymb,latexsym,epsfig,psfrag,graphicx}
\usepackage{tikz}
\usetikzlibrary{positioning,calc} 
\usepackage{pgfplots}
\pgfplotsset{compat=1.18}
\newcounter{definition}

\usepackage{mathtools}
\usepackage[inline]{enumitem}
\usepackage{nameref}
\usetikzlibrary{intersections,trees,arrows.meta}
\usepackage{accents}

\newcommand{\averageE}{\mathbb{E}}
\newcommand{\realR}{\mathbb{R}}

\DeclareMathOperator{\trace}{trace}

\makeatletter
\newcommand{\sbsection}[1]{%
  \def\@currentlabelname{#1}
  \section{#1}
}
\makeatother

\definecolor{blockcoloryellow}{cmyk}{0.03,0.03,0.12,0.0}
\definecolor{blockcolor}{cmyk}{0.12,0.04,0.08,0.0}
\definecolor{arred}{cmyk}{0,1.0,1.0,0.30}
\definecolor{arblue}{cmyk}{1.0,0.4,0,0}

\jname{To be published in Annual Review of Control, Robotics, and Autonomous Systems}
\jvol{10}
\jyear{2027}

\begin{document}

\markboth{Meijer et al.}{Dual Control}

\title{Dual~Control: On Exploration--Exploitation in~Linear~Systems}

\author{Tomas J. Meijer$^1$ and Anders Rantzer$^1$
\affil{$^1$Department of Automatic Control, Faculty of Engineering (LTH), Lund University, SE-221 00 Lund, Sweden; email: \{tomas.meijer,anders.rantzer\}@control.lth.se}}

\begin{abstract}
The term ``dual control'' refers to the dual objective of simultaneously balancing exploration and exploitation. Problems of this kind have been studied for nearly a century. This paper is devoted to theory and methodology relevant for optimal control of linear time-invariant systems whose parameters are initially unknown and must be learned by active probing. 
We review the main ideas underlying four major research directions: Multi-armed bandits, self-tuning regulators, regret rate minimizing controllers, and minimax optimal dual controllers. The first three have a long history and rich literature, whereas the fourth provides a promising framework for robust dual control.
\end{abstract}

\begin{keywords}
dual control, online learning, active learning, adaptive control, reinforcement learning, data-driven control
\end{keywords}
\maketitle


\section{THE EXPLORATION--EXPLOITATION TRADEOFF}
The interplay between exploration and exploitation has been studied and rediscovered repeatedly in different scientific fields. A good example is behavioral ecology, where scientists since the 1940-50s have studied the question: \emph{Should an animal that has found food continue searching locally or move elsewhere}? This question was initially studied by ecologists in a descriptive manner, but in the 1960s it was rephrased as an optimization problem~\cite{macarthur1966optimal} that could be addressed using mathematical tools~\cite{charnov1976optimal}. 
Other examples can be found in mining,~\cite{jung2021systematic}, neuroscience~\cite{cohen2007should} and organizational science~\cite{lavie2010exploration}. A broad view of the subject was provided in the monograph by Holland~\cite{holland1992adaptation}, first published in 1975.

Already in 1933, Thompson~\cite{thompson1933likelihood} gave a mathematical formulation of the multi-armed bandit problem in the context of clinical trials and introduced the idea now known as Thompson sampling. During World War II, Allied scientists studied the bandit problem for resource allocation and defense logistics. According to Whittle~\cite{whittle1979discussion}, \emph{``the problem is a classic one; it was formulated during the war, and efforts to solve it so sapped the energies and minds of Allied analysts that the suggestion was made that the problem be dropped over Germany, as the ultimate instrument of intellectual sabotage.''} This quote is from a publication devoted to Gittins' Bayesian approach to the bandit problem, which will be discussed in a later section.

\subsection{Feldbaum's Pioneering Work}
\begin{figure}[!bt]
\centering
\resizebox{.7\linewidth}{!}{%
\begin{tikzpicture}[
    block/.style={draw, thick, rectangle, minimum height=1.1cm, minimum width=2.4cm, align=center, fill=blockcolor,inner sep=6pt},
    sum/.style={draw, thick, circle, minimum size=0.7cm},
    arrow/.style={->, thick},
    node distance=.5cm and 1.8cm
]
    \node[block] (plant) {Uncertain linear system\\$x(t+1)=Ax(t) + Bu(t) + w(t)$\\$(A,B)\in\mathcal{M}$};
    \node[name=disturbance,left=of plant.170] {};
    \node[block, below=of plant] (controller) {Causal dual controller\\ $u(t)=\mu_t(x(0),\hdots,x(t),u(0),\hdots,u(t-1))$};
    \node[name=output, right=of plant.0] {};
    \coordinate[name=input, left=of plant.190];
    \coordinate[name=output2, right=of plant.7];
    \coordinate[name=input2, left=of plant.170];

    \draw[very thick,->] (disturbance) -- node[above,pos=0.1] {$w(t)$} (plant.170);
    \draw[very thick,->] (plant.0) -- node[above,pos=0.9] {$x(t)$} (output);
    \draw[very thick,->] ($(plant.0)!.7!(output)$) |- (controller.0);
    \draw[very thick,->] (controller.180) -| ($(input)!.3!(plant.190)$) -- node[above,pos=0.1] {$u(t)$} (plant.190);
\end{tikzpicture}%
}
\caption{In the dual control setting, the parameters $(A,B)\in\mathcal{M}$ are initially unknown and must be learned by the controller based on measurements of the state $x(t)$ while simultaneously regulating the plant through the control input $u(t)$. The different approaches covered in this review are distinguished by different assumptions on the admissible model class $\mathcal{M}$ and the exogenous disturbance $w(t)$.}
\label{fig:schematic}
\end{figure}
The first researcher to formulate a mathematical problem treating the exploration--exploitation tradeoff in its full generality was Feldbaum. He introduced the
term \emph{dual control} in the early 1960s~\cite{feldbaum1960dual,feldbaum1963dual} to describe the need for feedback controllers, such as the one shown in Figure~\ref{fig:schematic}, to simultaneously learn and regulate the plant.
In the following years, this idea propagated into a wide variety of subject areas in engineering, including adaptive control, reinforcement learning, and Bayesian optimization.

\begin{marginnote}[]
\entry{Dual control}{Feedback control that simultaneously regulates the plant and probes it to reduce uncertainty.}
\end{marginnote}
Feldbaum emphasized that learning often needs to be active: Without probing, you will not learn how the system responds. This makes dual control very different from classical communication theory and signal processing. Aggressive probing may deteriorate performance but improve parameter estimates. Conversely, cautious control may give better short-term performance but fail to gather information. Feldbaum addressed the dual control problem using Bellman's work
on dynamic programming~\cite{Bellman57}, by combining the physical state with an information state. He also recognized that solving Bellman's equation exactly poses computational challenges and emphasized the need for tractable approximations.

A few years later, in 1965, {\AA}str{\"o}m~\cite{aastrom1965optimal} introduced the notion of Partially Observable Markov Decision Processes (POMDPs) as a rigorous mathematical framework for the study of dual control. Such models have a discrete state space, which simplifies representation of the information state, and have later gained widespread recognition in computer science. However, the computational challenges persisted.

\begin{marginnote}[]
\entry{POMDP}{Markov decision process with incomplete state information. Decisions are based on a belief state.}
\end{marginnote}

\begin{figure}[!bt]
\centering
\resizebox{.8\linewidth}{!}{%
\begin{tikzpicture}[    
    grow                    = down,
    level distance          = 12mm,
    edge from parent path = {(\tikzparentnode.south) -- ++(0,-2mm) -| (\tikzchildnode.north)},
    edge from parent/.style = {draw, -},
    every node/.style       = {draw, align=center, font=\small},
    level 1/.style          = {sibling distance=62mm},
    level 2/.style           = {sibling distance=40mm},
    level 3/.style          =   {sibling distance=27mm},
    root/.style             = {draw=none, fill=none, font=\bfseries, text=arred},
    branch/.style           = {draw=none, fill=blockcoloryellow, text width=26mm, text=arred},
    leaf/.style              = {draw=none, fill=blockcolor, text width=20mm, font=\scriptsize, text=arblue},
    plainleaf/.style = {draw=none, fill=none, text width=20mm, font=\scriptsize},
]
    \node[root] {Dual Control in Linear Systems}
        child { node[branch] {Static plant}
            child[level distance=24mm, edge from parent path={
    (\tikzparentnode.south) -- ++(0,-18mm) coordinate (auxskip)
    (auxskip) -| (\tikzchildnode.north)
}] { node[leaf] {Multi-armed bandits}
                child[level distance=12mm, edge from parent/.style={draw, dotted}, edge from parent path={(\tikzparentnode.south) -- ++(0,-2mm) -| (\tikzchildnode.north)}]
      { node[plainleaf] {Section~\ref{sec:bandits}}}
            }
        }
        child { node [branch] {Dynamic plant} 
            child { node[branch] {Stochastic disturbances} 
                child { node[leaf] {Self-tuning regulator}
                    child[edge from parent/.style={draw, dotted}]  { node[plainleaf] {Section~\ref{sec:STR}}}
                }
                child { node[leaf] {Regret rate minimization}
                    child[edge from parent/.style={draw, dotted}]  { node[plainleaf] {Section~\ref{sec:regretrate}}}
                }
            }
            child { node[branch] {Worst case disturbances} 
                child { node[leaf] {Minimax dynamic games}
                    child[edge from parent/.style={draw, dotted}]  { node[plainleaf] {Section~\ref{sec:minimax}}}
                }
            }
        };    
\end{tikzpicture}%
}
\caption{Outline of the four main research directions considered in this review. Multi-armed bandits, while static, have contributed significantly to our understanding of the exploration--exploitation tradeoff. For dynamic systems with stochastic disturbances, we review the classical self-tuning regulator and more recent regret rate minimizing controllers. While the self-tuning regulator allows for colored noise, the regret rate minimizing controllers typically assume white noise. The last approach deals with dynamic systems subject to worst-case disturbances including colored noise.}
\label{fig:tree}
\end{figure}

\subsection{Further literature}
The focus of this review is on exploration--exploitation in linear systems, i.e., systems with a continuous state space. The distinction between discrete and continuous is somewhat ambiguous, however, for example due to the fact that probability distributions for discrete states have continuous parametrizations. 
Historically, literature on dual control of linear systems is closely tied to adaptive control, which is summarized in numerous surveys~\cite{kumar1985survey,wittenmark1995adaptive,annaswamy2021historical} and textbooks~\cite{aastrom2013adaptive,goodwin2014adaptive,sastry2011adaptive,chen2012identification,astolfi2007nonlinear}. In particular, adaptive control has triggered several attempts to derive exact solutions for special cases of dual control. Examples include~\cite{sternby1976simple,86aastrom+,bernhardsson1989dual} and~\cite{kulcsar1996dual}. In this review, the treatment of adaptive control will be restricted to self-tuning regulators, as a valuable benchmark for comparison with dual control algorithms.

In recent years, dual control and the exploration--exploitation tradeoff have mainly been discussed in the context of reinforcement learning and 
data-driven control. Most of the reinforcement learning literature is developed in a discrete-time setting, see, e.g., Sutton \& Barto~\cite{sutton1998reinforcement}, but there is also a rich literature studying linear systems from the perspective of statistical machine learning. See, for instance, the survey paper~\cite{tsiamis2023statistical}. In parallel, dual model predictive control (MPC) has also become an active research area, which will not be the focus of this review but was reviewed by Mesbah~\cite{mesbah2018stochastic} in 2018. Dual MPC aims to approximate the exploration--exploitation tradeoff through receding-horizon optimization. These approaches have the benefit that, unlike most dual controllers, they are also able to handle actuator and performance/safety constraints. Some noteworthy dual MPC approaches include approximate dynamic programming methods, explicit dual-control formulations, and stochastic output-feedback MPC schemes that seek to preserve and exploit the dual control effect~\cite{Parsi2023a,Parsi2023b,Arcari2020,Messerer2022}. More recent developments focus on, e.g., leveraging posterior Gaussian process covariance for active learning, quasi-linear parameter-varying systems, and promoting exploration by incorporating a covariance-dependent quadratic surrogate for the information gain into the MPC stage cost.~\cite{Baltussen2025b,Baltussen2026,Mulagaleti2026}.

Given this background, this review focuses, as summarized in Figure~\ref{fig:tree}, on the following four main research directions:
\begin{enumerate*}[label=(\alph*)]
    \item Multi-armed bandits,
    \item self-tuning regulators,
    \item regret rate minimization for linear-quadratic control, and 
    \item minimax optimal dual control.
\end{enumerate*}
The first three are included for their historical importance to dual control theory, while the last direction, rooted in the theory of \emph{zero sum dynamic games}, is presented as an exciting direction for ongoing and future research.

\section{THEORY FOR MULTI-ARMED BANDITS}\label{sec:bandits}
The most well-developed theoretical framework for understanding of the exploration--exploitation tradeoff is the theory for multi-armed bandits. The problem can be formulated as follows: Consider a slot machine (bandit) with $N$ arms. For each $i\in\{1,\ldots,N\}$, the arm $i$, if pulled, has a fixed probability $\theta_i$ of yielding a reward of one unit. Otherwise, the reward is zero. The parameter vector $\theta=(\theta_1,\ldots,\theta_N)$ is unknown, but in each discrete time instance $t$ one arm $i(t)$ should be pulled, resulting in the reward $x_{i(t)}(t)$. The goal is to maximize the expected reward over time. The problem is static in the sense that all arms are memoryless, but the exploration--exploitation tradeoff appears in the dynamic interaction between information and rewards. 

\begin{marginnote}[]
\entry{Multi-armed bandit}{Sequential decision-making problem balancing exploration and exploitation for maximal accumulated reward.}
\end{marginnote}

\subsection{A Bayesian setting: The Gittins index}
A sequence of research contributions initiated by Gittins \& Jones~\cite{gittins1974dynamic} solved the multi-armed bandit problem in a Bayesian setting. A prior distribution for each $\theta_i$ is supposed to be given. The goal is here to maximize the total expected reward $\averageE\sum_{t=0}^\infty\beta^tx_{i(t)}(t)$, where $\beta\in(0,1)$ is a discount factor. Let $q^i(t)$ denote the conditional probability density of $\theta_i$, given outcomes up to time $t$. A function $\nu$, called the \emph{Gittins index}, can be defined on the set of probability densities such that an optimal strategy is obtained by always pulling the arm for which $\nu(q^i(t))$ takes the largest value. 

\begin{marginnote}[]
\entry{Gittins index}{Bayesian priority index that yields the optimal discounted policy for classical multi-armed bandits.}
\end{marginnote}

An intriguing aspect of the solution by Gittins \& Jones is its decentralized nature. Conditional probability densities can be updated for each arm separately and then used for evaluation of the Gittins index. Finally, the values are compared to decide which arm to activate next. Without going into details, the function $\nu$ can be thought of as the maximum expected reward per unit of time, under the assumption that the game can be stopped anytime. Thus, computation of the Gittins index is an \emph{optimal stopping} problem.

The exploration--exploitation nature of optimal solutions to the multi-armed bandit problem is further highlighted by a result due to Kelly~\cite{kelly1981multi}: As the discount factor $\beta$ approaches one, the optimal policy tends to initially follow the \emph{rule of least failures}, i.e., use the arm which has incurred the least number of failures (and in case of non-uniqueness the one among the machines with fewest failures one with most successes). This rule naturally leads to \emph{exploration}. However, the optimal policy will eventually stop switching and play only one arm, hence practicing \emph{exploitation}.

Finally, it should be noted that application of the Gittins index is not restricted to multi-armed bandits. Instead of defining it as a function of conditional probability densities for bandits, it could be defined on the state of any Markov process. Then the theory applies as long as there are $N$ copies of the process and at every time instant the permissible control action is to let one of the processes evolve while freezing all the rest.

\subsection{Regret optimization for multi-armed bandits}

A decade later, Lai \& Robbins~\cite{lai1985asymptotically} addressed the multi-armed bandit using a non-Bayesian perspective, with the objective of minimizing the expected difference between bandit reward for the optimal arm and the reward of a policy without prior knowledge of $\theta$. This difference
\begin{align*} 
    R_T:=\underbrace{\max_i\averageE\sum_{t=0}^Tx_{i}(t)}_{\text{\begin{tabular}{l}Expected reward using optimal arm\end{tabular}}}-\underbrace{\averageE\sum_{t=0}^Tx_{i(t)}(t)}_{\text{\begin{tabular}{l}Expected reward using policy\end{tabular}}}
\end{align*}
is called the \emph{regret} of the policy and is usually studied with focus on large values of $T$. Lai \& Robbins~\cite{lai1985asymptotically} proved that $R_T$ must grow asymptotically at least with the rate $\log T$ and they also constructed a policy achieving this lower bound asymptotically. The policy of Lai \& Robbins has a form similar to the Bayesian setting: A function called the \emph{upper confidence index} or upper confidence bound is evaluated for each arm and the policy is to pick the arm with the highest value of the index.

\begin{marginnote}[]
\entry{Regret}{The performance loss relative to the optimal policy that knows the true parameters.}
\end{marginnote}

\subsection{Optimism in the face of uncertainty for bandits}
The upper confidence index of Lai \& Robbins was later simplified by Agrawal~\cite{agrawal1995sample} to 
\begin{align*}
    \hat\theta_i(t)+\sqrt{(2\log t)/n_i(t)}.
\end{align*}
Here, $n_i(t)$ is the number of times arm $i$ has been pulled by time $t$ and $\hat\theta_i(t)$ is the corresponding empirical mean reward.
The policy was also proved by Auer, \mbox{Cesa-Bianchi} \& Fischer~\cite{auer2002finite} to have logarithmic regret uniformly, not only asymptotically. The second term in the index quantifies the uncertainty in the empirical mean.
Choosing the arm that maximizes the upper confidence bound then leads to a policy that exhibits ``\emph{optimism in the face of uncertainty}'', promoting exploration at the expense of exploitation.

\begin{marginnote}[]
\entry{Optimism in the face of uncertainty}{When different decisions are compared, give preference to those that have been tried less in the past.}
\end{marginnote}

\subsection{Minimax optimal control of bandits}
Also bandits have been studied from the minimax optimal control perspective. In~\cite{auer1995gambling}, there are \emph{no statistical assumptions} on the rewards. Instead regret is defined as 
\begin{align*} 
    R_T:=\max_j\sum_{t=0}^Tx_{j}(t)-\sum_{t=0}^Tx_{i(t)}(t)\
\end{align*}
where an adversarial player selects a reward vector $(x_1(t),\ldots,x_N(t))\in[0,1]^N$ arbitrarily. 
With a deterministic policy, not much can be done. 
However, with a randomized policy, $R_T$ becomes a random variable and minimization of $\averageE R_T$ leads to an interesting exploration--exploitation tradeoff \cite{auer2002finite}.

\begin{textbox}[!bt]
	\sbsection{SELF-TUNING CONTROL}
	\label{sb:str}
	We introduce the idea behind self-tuning regulators in a simple setting. Consider the input-output model 
	\begin{align*}
		y(t+1)&=-a_0y(t)-\ldots -a_{n-1}y(t-n+1)
		+b_0u(t)+\ldots +b_{n-1}u(t-n+1)+w(t)
	\end{align*}
	where $u$ is the input, $y$ is the output, and $w$ is a disturbance. If $w$ is zero-mean white noise and the system parameters are known, a \emph{minimum-variance controller} is given by 
	$$
	u(t)=\frac{1}{b_0}\left(\sum_{k=0}^{n-1}a_ky(t-k)-\sum_{k=1}^{n-1}b_ku(t-k)\right).
	$$
	This policy achieves $y=w$, so that $\averageE y^2(t)=\averageE w^2(t)$ which is the best we can hope for when $w$ is unmeasurable white noise.
	The self-tuning regulator uses the same structure, but replaces the system parameters by estimates that are updated recursively. 
	Specifically, assuming that $b_0$ is known, the parameter estimate 
	$
	\hat{\theta}(t)\coloneqq(\hat{a}_0(t),\ldots,\hat{a}_{n-1}(t),\hat{b}_1(t),\ldots,\hat{b}_{n-1}(t))
	$
	is obtained as the minimizer of the quadratic expression
	\begin{eqnarray*}
		\sum_{\tau=n-1}^{t-1}\left(y(\tau+1)+\sum_{k=0}^{n-1}a_ky(\tau-k)-\sum_{k=0}^{n-1}b_ku(\tau-k)\right)^2.
	\end{eqnarray*}
	For Gaussian white noise $w$, $\hat{\theta}(t)$ is the maximum likelihood estimate. Let $\phi(t)\coloneqq (-y(t),\ldots,-y(t-n+1),u(t-1),\ldots,u(t-n+1))$. The self-tuning controller is defined recursively as $u(t)=-b_0^{-1}\hat{\theta}^\top(t)\phi(t)$ with 
	\begin{eqnarray*}
		P(t)&=&P(t-1)-\frac{P(t-1)\phi(t)\phi^\top(t) P(t-1)}{1+\phi^\top(t) P(t-1)\phi(t)},\\
		\hat{\theta}(t)&=&\hat{\theta}(t-1)+\frac{P(t-1)\phi(t-1)}{1+\phi^\top(t)P(t-1)\phi(t)}\left[y(t)-\phi^\top(t-1)\hat{\theta}(t-1)\right],
	\end{eqnarray*}
	with initial conditions such that $P(t)=\left(\sum_{\tau=n-1}^{t}\phi(\tau)\phi^\top(\tau)\right)^{-1}$ is well defined. 
	This algorithm gives $\theta(t)-\hat\theta(t)$ covariance $P(t)$ when $w$ is white noise with unit variance.
	The original algorithm of~\cite{ast+wit73} was actually more general, with $w$ not necessarily white noise. Moreover, a fixed number of $b$-coefficients were allowed to be zero, corresponding to a pre-specified delay from inputs to outputs in the model. 
	Such self-tuning controllers were used successfully in various applications as well as commercial products. See, for example, the book by {\AA}str{\"o}m \& Wittenmark~\cite[Chapter~12]{aastrom2013adaptive} or the survey~\cite[Section 5]{annaswamy2021historical}.
\end{textbox}

\section{SELF-TUNING REGULATORS FOR LINEAR SYSTEMS}\label{sec:STR}
We will now shift the discussion to control of linear time-invariant \emph{dynamical} systems. It is natural to start with \emph{self-tuning regulators}. This concept, introduced by {\AA}strom \& Wittenmark~\cite{ast+wit73}, became a dominant paradigm of \emph{adaptive control} during the 1970s-1980s. The algorithm was also implemented in several commercial products. See \cite[Chapter~12]{aastrom2013adaptive}.

As detailed in the textbox~``\nameref{sb:str}'', the main idea of the self-tuning regulator was to use the so-called minimum-variance control setting, but replace the system parameters by estimates that are updated recursively based on past data. The focus was not to optimize the exploration--exploitation tradeoff, but several aspects related to dual control naturally emerged in the analysis.

\begin{marginnote}[]
\entry{Self-tuning regulator}{Minimum variance controller that recursively estimates model parameters and redesigns the controller on-line.}
\end{marginnote}


\subsection{Design considerations and fundamental limitations}
In this section, we discuss some important issues and insights related to dual control that naturally emerged in the analysis of self-tuning regulators.
\subsubsection{The $B$-polynomial and the self-tuning control objective}
In our description of the algorithm, the parameter $b_0$ was supposed to be known. In practical applications, however, no system follows the mathematical model exactly, so $b_0$ and a potential input-output delay should not necessarily be considered as prior knowledge about the process. Instead, they are tuning parameters for the controller. Specifically, $b_0$ quantifies the gain of the controller, i.e., how aggressive it is. 
Similarly, time delays in the model quantify how fast disturbances are rejected. 

Even when the model structure perfectly matches the real process, the control objective could lead to unexpected consequences. The minimum-variance controller is a ``dead beat'' controller, in the sense that it tries to cancel disturbances as quickly as possible. If the sampling rate is high, this typically leads to overly-aggressive control action. Hence, the sampling time is an important design parameter in the context of self-tuning control.

Yet another perspective on minimum-variance control is that it specifies the closed-loop transfer function without consideration of the open-loop plant characteristics. This leads to pole-zero cancellations and the resulting closed-loop system can become internally unstable. For example, consider the following simple case with a non-minimum-phase zero:
\begin{align*}
    y(t+1)&=y(t)-u(t)+2u(t-1)+w(t).
\end{align*}
The minimum-variance controller $u(t)=y(t)+2u(t-1)$ gives $y(t+1)=w(t)$, so the input grows exponentially as $u(t+1)=2u(t)+w(t)$. To prevent this, a standard assumption is that all roots of the polynomial $b_0z^{n-1}+\cdots+b_{n-2}z+b_{n-1}$ must be inside the unit disc.

\subsubsection{The $C$-polynomial and the effect of colored noise}
An important motivating application for the development of self-tuning control was disturbance attenuation in process control. The disturbances may enter the process in many different ways and the assumption that $w(t)$ is white noise can lead to biased estimates and potentially poor control. To avoid such effects, it is common to use the ARMAX model
\begin{marginnote}[]
\entry{ARMAX model}{AutoRegressive Moving-Average model with exogenous inputs.}
\end{marginnote}
\begin{eqnarray*}
    y(t+1)&=&-\sum_{k=0}^{n-1}a_ky(t-k)
    +\sum_{k=0}^{n-1}b_ku(t-k)+e(t+1)+\sum_{k=0}^{n-1}c_ke(t-k)
\end{eqnarray*}
where $e$ is assumed to be white noise. This model is obviously much more flexible than the previous one, since the white noise $w$ now is replaced by a moving average of $e$. Assume that a minimum-variance controller for the ARMAX model removes the effect of all disturbances except the last one, i.e., $y(t)=e(t)$ for all $t$. Then, the closed-loop dynamics becomes
\begin{eqnarray*}
    y(t+1)&=&-\sum_{k=0}^{n-1}(a_k+c_k)y(t-k)
    +\sum_{k=0}^{n-1}b_ku(t-k)+e(t+1),
\end{eqnarray*}
where the bias effect in parameter estimation is reflected by the fact that $a_k$ is replaced by $a_k+c_k$. Connecting the corresponding controller 
\begin{eqnarray*}
    u(t)=\frac{1}{b_0}\left(\sum_{k=0}^{n-1}(a_k+c_k)y(t-k)-\sum_{k=1}^{n-1}b_ku(t-k)\right)
\end{eqnarray*}
to the ARMAX model yields 
\begin{eqnarray*}
    y(t+1)&=&-\sum_{k=0}^{n-1}c_ky(t-k)+e(t+1)+\sum_{k=0}^{n-1}c_ke(t-k),
\end{eqnarray*}
which is stable only if all roots of the polynomial $z^n+c_0z^{n-1}+\cdots+c_{n-2}z+c_{n-1}$ are inside the unit disc. This illustrates that the assumption about removing ``the effect of all disturbances except the last one'' may not be achievable by a stabilizing controller and extra assumptions are needed for convergence of the self-tuning controller. Following relevant counterexamples in~\cite{ljung1975counterexamples}, Ljung~\cite{ljung77on} formulated the stronger condition that
\begin{equation}\label{eq:Cinv_min_half}
    (1+c_0z^{-1}+\cdots+c_{n-2}z^{-(n-1)}+c_{n-1}z^{-n})^{-1}-1/2
\end{equation}
should be positive real, which is now standard in the literature on self-tuning control.

\subsection{Persistence of excitation and robustness to unmodelled dynamics}
A central component of the self-tuning controller is parameter estimation using least-squares optimization. Accurate estimation is only possible if the data is sufficiently rich. This is true in all areas of statistics, but some aspects need extra attention when data is generated in closed loop. The reason is once again the tradeoff between exploration and exploitation: If only the optimal controller is used, then we do not collect information about other options.

In 1966, {\AA}str{\"o}m \& Bohlin~\cite{aastrom1965numerical}
introduced the following definition of persistence of excitation. A bounded signal $u$ is said to be \emph{persistently
    exciting of order $m$} if the limits
    \begin{eqnarray*}
        \lim_{N\to\infty}\frac{1}{N}\sum_{t=1}^Nu(t),\text{ and }
        r_u(T):=\lim_{N\to\infty}\frac{1}{N}\sum_{t=1}^Nu(t)u(t+T),
    \end{eqnarray*}
    exist and the matrix $R_u:=\{r_u(i-j) : i,j=1,\ldots,m+1\}$ is positive definite.

\begin{marginnote}[]
\entry{Persistence of excitation}{Input ``richness'' condition that, combined with controllability assumptions, ensures that maximum likelihood estimates of ARMAX coefficients are consistent.}
\end{marginnote}
Other versions of the persistence of excitation condition were later used to prove exponential stability of adaptive control schemes, for example in continuous time by Anderson~\cite{anderson1977exponential} and Morgan \& Narendra~\cite{morgan1977stability} and discrete time by Bai \& Sastry~\cite{bai1985persistency}. Exponential stability is a very valuable property in systems analysis, since it is a natural starting point for the study of robustness to unmodelled dynamics. However, provable excitation levels and exponential convergence rates are often established using conservative bounds, which in the dual control context could lead to overly-explorative policies at the expense of exploitation.

\section{REGRET RATE MINIMIZATION FOR LINEAR-QUADRATIC CONTROL}\label{sec:regretrate}
The modern literature on feedback control with on-line learning is to large extent inspired by developments in reinforcement learning, where statistical concentration bounds and asymptotic regret rates play a dominant role. See Recht~\cite{recht2019tour} or Tsiamis et al.~\cite{tsiamis2023statistical} for overviews.

\subsection{Logarithmic regret bounds for self-tuning control}\label{sec:regret-bounds-str}
Regret analysis of adaptive schemes can be traced back to the mid 1980s, when Lai \& Wei developed theory for adaptive controllers in parallel with multi-armed bandits. As illustrated in the textbox~ ``\nameref{sb:log-regret}'', in many cases regret has a logarithmic growth rate. Specifically, it turns out that logarithmic growth rate of $R_T$ generally holds in the context of self-tuning control. Results of this kind were pioneered by Lai \& Wei~\cite{lai1987asymptotically} and further refined by Guo~\cite{guo1995convergence}. In 2020, Cassel et al.~\cite{cassel2020logarithmic} proved logarithmic regret growth also for more general linear-quadratic control problems with a state-space model $(A,B)$ where the $B$-matrix is known and $A$ can be estimated from intrinsic noise. From a dual control perspective the cases where significant exploration is needed to find the optimal controller are more interesting. In such cases, logarithmic growth cannot be achieved and the regret instead grows as $\sqrt{T}$. The reasons for this will be elaborated next.

\begin{textbox}[!bt]
\sbsection{WHY LOGARITHMIC REGRET APPEARS}
\label{sb:log-regret}
    The following simple example by Lai \& Robbins~\cite{lai1978adaptive} explains how logarithmic growth appears naturally. Consider the (static) linear model
    \begin{eqnarray*}
        y(t)&=&bu(t)-\mu+w(t),
    \end{eqnarray*}
    where $b\ne0$ is known, while $\mu$ is an unknown parameter and $w(t)$ are independent identically distributed random variables with zero mean and unit variance. The input $u(t)$ is used in order to minimize $\averageE y(t)^2$. After collecting values of $u(k)$ and $y(k)$ for $k=1,\ldots,t$, least-squares estimation of $\mu$ gives
    \begin{eqnarray*}
        \hat\mu(t):=\frac{1}{t}\sum_{k=1}^t[bu(k)-y(k)]
        =\mu-\frac{1}{t}\sum_{k=1}^tw(k),\text{ and } 
        \averageE[\hat\mu(t)-\mu]^2=\averageE\left[\frac{1}{t}\sum_{k=1}^tw(k)\right]^2=\frac{1}{t^2}\sum_{k=1}^t\averageE w(k)^2=\frac{1}{t}.
    \end{eqnarray*}
    Minimizing the objective $\averageE y(t)^2$ yields the control law $u(t)=\hat\mu(t-1)/b$ for $t\geqslant2$, which gives
    \begin{eqnarray*}
        \underbrace{\averageE \sum_{t=1}^Ty(t)^2}_{\text{\begin{tabular}{l}Expected cost\\using policy\end{tabular}}}-\underbrace{\averageE \sum_{t=1}^Tw(t)^2}_{\text{\begin{tabular}{l}Optimal\\expected cost\end{tabular}}}&=&
        \averageE \sum_{t=1}^T[bu(t)-\mu]^2
        =\averageE[bu(1)-\mu]^2+\averageE\sum_{t=2}^T[\hat\mu(t-1)-\mu]^2\\[-8mm]
        &=&\averageE[bu(1)-\mu]^2+1+\frac{1}{2}+\frac{1}{3}+\cdots+\frac{1}{T-1},
    \end{eqnarray*}
    which grows as the \emph{logarithm of $T$}. As before, the left-hand side can be viewed as the \emph{regret} $R_T$, quantifying the cost difference with and without knowledge of the parameter value $\mu$.
\end{textbox} 
\subsection{Regret rate for general linear-quadratic control}
The first general learning algorithm to optimize a quadratic objective for general linear systems with a regret bound growing as $\sqrt{T}$ was proposed in 2011 by Abbasi-Yadkori \& Szepesvári~\cite{abbasi2011regret}. The algorithm modified earlier work from 1998 by Campi \& Kumar~\cite{campi1998adaptive} using the idea of ``optimism in the face of uncertainty'' from multi-armed bandits. For linear systems, the idea was to construct high-probability confidence sets around the model parameters, find the optimal controller for each member of the confidence set, and finally choose the controller whose associated average cost is the smallest. This stimulates exploration when parameters are uncertain. The algorithm by Abbasi-Yadkori \& Szepesvári was not computationally efficient, but inspired a long sequence of follow-up contributions competing to find algorithms with the best possible bounds on regret.

It took some time to uncover reasons why some problems can be solved with $\log T$ regret bounds, while others need $\sqrt{T}$. Eventually, in 2020, this was settled in three parallel contributions~\cite{cassel2020logarithmic,simchowitz2020naive,ziemann2020phase}. More recently, these results were also extended to partially-observed (output feedback) systems~\cite{ziemann2025regret}. The difference between the two problem classes can be understood as follows: It has long been recognized that system parameters can be hard to identify in closed-loop operation ~\cite{soderstrom1975identifiability}. In particular, optimization of a linear system with a quadratic objective leads to a closed loop where not all parameters are identifiable. If the non-identifiable parameters are unimportant, this issue can be ignored and $\log T$ regret rate can be achieved. Conversely, if the optimal controller makes essential parameters non-identifiable, which is often the case, then exploration needs to drive the system away from optimality with higher regret as a result.

The mechanism can be illustrated by returning to the model $y(t)=bu(t)-\mu+w(t)$ from~``\nameref{sb:log-regret}'', now with uncertainty in both $b$ and $\mu$. Data collected with a constant control value $u$ makes it impossible to distinguish the pair $(b,\mu)$ from the pair $(b+\delta,\mu+u\delta)$, where $\delta$ is any real number. If the objective is to minimize $\averageE y(t)^2$, then the optimal policy corresponding to $(b,\mu)$ is $u(t)=\mu/b$, while the optimal policy corresponding to $(b+\delta,\mu+u\delta)$ is $u(t)=(\mu+u\delta)/(b+\delta)$. However, $u_*=\mu/b$ gives $u_*=(\mu+u_*\delta)/(b+\delta)$, so non-identifiability of the true parameters under the optimal policy does not create a need for exploration away from optimality. This is consistent with both the logarithmic regret growth that was calculated in Section~\ref{sec:regret-bounds-str} and the logarithmic regret of self-tuning control, which has a similar structure.

If the objective is changed to minimization of $\averageE[y(t)^2+u(t)^2]$, however, then the optimal policy for $(b,\mu)$ becomes $u(t)=b\mu/(1+b^2)$, while for $(b+\delta,\mu+u\delta)$ it is $u(t)=(b+\delta)(\mu+u\delta)/(1+(b+\delta)^2)$. In this case, the non-identifiability matters and exploration away from optimality is needed. 
This is consistent with the fact that the best achievable growth rate for the regret in general linear-quadratic control problems is $\sqrt{T}$ rather than $\log T$. Finally, it should  be noted that these rate bounds rely on idealized conditions and even small mismatch between the model structure and the true system will typically lead to growth rate $T$ instead of $\sqrt{T}$~\cite{lee2024nonasymptotic}.

\subsection{Certainty-equivalence control with additional probing}
After extensive efforts with advanced exploration schemes inspired by the reinforcement learning literature, such as ``optimism in the face of uncertainty'' and ``Thompson sampling''~\cite{ouyang2017control}, it finally turned out that also the classical certainty-equivalence idea of Bar-Shalom \& Tse~\cite{bar1974dual}, dating back to 1974, is able to maintain regret growth rates matching the theoretical lower bounds~\cite{simchowitz2020naive}. 

\begin{marginnote}[]
\entry{Certainty-equivalence control}{Control strategy that uses current estimates as if they were the true parameters.}
\end{marginnote}
Specifically, Jedra \& Proutiere~\cite{jedra2022minimal} consider control of a linear system $x(t+1)=Ax(t)+Bu(t)+w(t)$ with the objective to minimize $\averageE\sum_{t=1}^T\left[x^\top(t) Qx(t)+u^\top(t) Ru(t)\right]$ using a causal adaptive control policy $\pi$. The system is supposed to be driven by zero-mean white noise $w$ with variance $\sigma^2I$. The matrices $A,B$ are unknown, but the policy has access to a controller $K_0$ that is stabilizing for all $A,B$. The policy of~\cite{jedra2022minimal} then has the following form: Given the data $x(0),u(0),\ldots,x(t-1),u(t-1),x(t)$, compute a least-squares estimate $(\hat{A}_t,\hat{B}_t)$ of $(A,B)$, then compute $K_t$ by solving the corresponding Riccati equation, and let
\[
    u(t):=\begin{cases}
        K_tx(t)+\nu(t),&\text{if past data is sufficiently rich},\\
        K_0x(t)+\nu(t),&\text{otherwise},
    \end{cases}
\]
where
\begin{eqnarray*}
    \begin{bmatrix}\hat{A}_t&\hat{B}_t\end{bmatrix}
    &:=&\left(\sum_{k=0}^{t-2}x(k)
    \begin{bmatrix}x(k)\\u(k)\end{bmatrix}^\top\right)
    \left(\sum_{k=0}^{t-2}\begin{bmatrix}x(k)\\u(k)\end{bmatrix}
    \begin{bmatrix}x(k)\\u(k)\end{bmatrix}^\top\right)^{-1},\\
    K_t&:=&-(R+\hat{B}_t^\top \hat{P}_t\hat{B}_t)^{-1}\hat{B}_t^\top \hat{P}_t\hat{A}_t\;\text{ with }    \hat{P}_t=Q+\hat{A}_t^\top\left(\hat{P}_t^{-1}+\hat{B}_tR^{-1}\hat{B}_t^\top\right)^{-1}\hat{A}_t,
\end{eqnarray*}
Excitation is provided by a sequence of independent random vectors $\nu(t)\in\mathcal{N}(0,\sigma^2\sqrt{d_x/t})$ where $d_x$ is the state dimension. Jedra \& Proutiere proved that with a proper definition of ``sufficiently rich'' data, the expected number of times that $K_0$ is used is \emph{finite}, while the expected value of the regret $R_T$ satisfies a bound of the form
\begin{eqnarray*}
    \averageE R_T&\leqslant C_1\sqrt{d_x}(d_x+d_u)\sqrt{T}\log T+C_2.
\end{eqnarray*}
Here $C_1$ and $C_2$ depend on the (unknown) values of $A$ and $B$. The resulting policy is a certainty-equivalence controller.
Interestingly, the upper bound on $\averageE R_T$ matches the lower bounds discussed before, so no other policy can do better in terms of expected asymptotic growth rate for general linear-quadratic problems. However, the specification of $\nu$ excludes logarithmic growth rate when the policy is applied to simpler problems such as self-tuning control, so in this sense the policy is not optimal.




Finally, let us note that the assumption that a stabilizing controller $K_0$ is known in advance is more restrictive than it might seem at first glance: Consider a scalar system $x(t+1)=ax(t)+bu(t)$ with unknown parameters $a\in[a_1,a_2]$ and $b\in[b_1,b_2]$. If $k_0$ is a stabilizing feedback, then the pair $(a,b)$ must be restricted to a set between two parallel lines in the parameter space, i.e., $-1\leqslant a+bk_0\leqslant 1$. In particular, if $b_1\leqslant 0\leqslant b_2$, then $-1<a_1\leqslant a_2<1$ and, thus, the system must be stable even without feedback.
This is in contrast with the next section, where unstable systems will be treated without prior knowledge of the sign of the $b$-parameter.


\section{MINIMAX OPTIMAL DUAL CONTROL FOR LINEAR SYSTEMS}\label{sec:minimax}
The final approach to dual control that we highlight in this review is based on minimax optimal control. While this approach is less mature than the ones discussed so far, it is a promising and active research area with many open issues that remain to be investigated. The minimax approach accounts for the worst-case of both disturbances and uncertain parameters, as discussed in~\cite{Cusumano/P88,Ioannou1988,vinnicombe2004examples,2004megretskinonlinear}. The minimax adaptive/dual control paradigm was first introduced for linear systems by Didinsky \& Basar in 1994~\cite{Didinsky1994} and for nonlinear systems in~\cite{Pan1998}. This setting can be viewed as a natural extension of the well-known $\mathcal{H}_{\infty}$ control problem~\cite{Zames1981,zhou+96} to an adaptive setting. 

\begin{marginnote}[]
\entry{Minimax control}{Optimize performance against worst-case disturbances and, in the context of adaptive/dual control, the worst-admissible model.}
\end{marginnote}

Consider an uncertain discrete-time linear time-invariant system described by
\begin{equation}\label{eq:system}
x(t+1) = Ax(t) + Bu(t) + w(t),\qquad x(0)=x_0,
\end{equation}
with state $x(t)\in\mathbb{R}^{n}$, input $u(t)\in\mathbb{R}^{m}$ and disturbance $w(t)\in\mathbb{R}^{n}$ at time $t\in\mathbb{N}$. The system parameters $A\in\mathbb{R}^{n\times n}$ and $B\in\mathbb{R}^{n\times m}$ are unknown but they belong to some known model class $\mathcal{M}\subset\mathbb{R}^{n\times n}\times \mathbb{R}^{n\times m}$. In the adaptive minimax optimal control problem considered here, we are interested in computing the worst-case quadratic cost across all $(A,B)\in\mathcal{M}$ and external disturbances $w$: For positive definite matrices $Q$ and $R$ and $\gamma>0$, we solve
\begin{equation}\label{eq:minimax-cost}
    J^\star(x_0)\coloneqq \inf_{\mu}\sup_{w,N}\sup_{(A,B)\in\mathcal{M}}\sum_{t=0}^{N-1} \|x(t)\|^2_Q + \|u(t)\|^2_R - \gamma^2\|w(t)\|^2,
\end{equation}
for the uncertain system in Equation~\eqref{eq:system} with $(A,B)\in\mathcal{M}$ and with a causal control policy 
\[
    u(t)=\mu_t(x(0),\hdots,x(t),u(0),\hdots,u(t-1)).
\]
The above problem can be viewed as a two-player dynamic game, in which one player, the controller, aims to minimize the cost while the other player chooses $w$, $A$, and $B$ to maximize the cost. This is the standard game formulation of the model-based $\mathcal{H}_{\infty}$ control problem~\cite{Basar/B95} except that the adversary also selects the unknown parameters $(A,B)\in\mathcal{M}$. Since the matrices are unknown but constant, the optimal feedback law is naturally incentivized to \emph{explore} so that it can ``learn'' $A$ and $B$ early on, in order to \emph{exploit} this knowledge later. Such nonlinear adaptive controllers are capable of stabilizing and optimizing the behavior even when no linear controller can simultaneously stabilize the system in Equation~\eqref{eq:system} for all $(A,B)\in\mathcal{M}$. Similar to $\mathcal{H}_{\infty}$ optimal control, $\gamma$ serves as a bound on the ``gain'' from disturbances $w$ to states $x$ and inputs $u$, i.e.
\begin{align*}
    \sum_{t=0}^{\infty} \|x(t)\|_Q^2 + \|u(t)\|^2_R 
    &\leqslant \gamma^2\sum_{t=0}^{\infty} \|w(t)\|^2&
    &\text{ for all }(A,B)\in\mathcal{M}.
\end{align*}
Although the minimax approach is known to be conservative when applied to systems with Gaussian white noise, even in the model-based setting, it offers several significant benefits compared to the other dual control approaches discussed before. As we will see, results for the minimax approach do not rely on the availability of a pre-stabilizing controller. This not only offers practical benefits but it also means that minimax dual control can deal with settings where no single static feedback controller can stabilize the system for all parameter values. 

\subsection{Fundamental limitations of minimax dual control}
Intuitively, it should be clear that no minimax dual controller can outperform the model-based $\mathcal{H}_{\infty}$ controller in terms of the cost in Equation~\eqref{eq:minimax-cost}, or achieve a better $\ell_2$-gain from disturbances to states than the gain achievable by a model-based $\mathcal{H}_\infty$ controller. Several works further explore fundamental limitations of the minimax dual control approach. In particular, Vinnicombe~\cite{vinnicombe2004examples} investigated the achievable closed-loop $\ell_2$-gain for scalar systems, i.e., $x(t+1) = ax(t) + bu(t) + w(t)$, where $x(t),u(t),w(t)\in\mathbb{R}$. Two distinct settings were considered: \begin{enumerate*}[label=(\roman*)]
    \item\label{item:uncertain-b} sign uncertainty in $b$, i.e., $\mathcal{M}=\left\{(a,1),(a,-1)\right\}$ with known $a$, and 
    \item\label{item:uncertain-a} sign uncertainty in $a$ with $b=1$, i.e., $\mathcal{M}=\left\{(a,1),(-a,1)\right\}$ and known $a$.
\end{enumerate*}

\begin{marginnote}[]
\entry{$\ell_2$-gain}{Worst-case amplification from disturbance energy to state/output energy.}
\end{marginnote}
\subsubsection{Sign uncertainty in the $B$-matrix}
For case~\ref{item:uncertain-b}, there exists no linear time-invariant controller that simultaneously stabilizes both systems in $\mathcal{M}$ for $a>1$. However, Vinnicombe~\cite{vinnicombe2004examples} showed that the nonlinear switching controller
\begin{equation}\label{eq:vinnicombe04controller} 
    u(t) = \begin{cases}
        -ax(t),\quad &\text{if }\sum_{\tau=0}^{t-1} (x(\tau+1)-ax(\tau))u(\tau)\geqslant 0,\\
        ax(t),\quad &\text{otherwise},
    \end{cases}
\end{equation}
stabilizes both systems and also achieves a closed-loop $\ell_2$-gain from $w$ to $x$ of 
$
\gamma^\star \coloneqq a + \sqrt{1+a^2}$. He also conjectured that no controller can achieve better $\ell_2$-gain for both systems, i.e., that the given controller is optimal. Recently, Rantzer~\cite{Rantzer2025acc} proved this conjecture to be true by showing that the cost for the infinite-horizon minimax optimal control problem is unbounded when $\gamma<\gamma^\star$.

\subsubsection{Sign uncertainty in the $A$-matrix}
For case~\ref{item:uncertain-a}, Vinnicombe~\cite{vinnicombe2004examples} showed that the adaptive controller in Equation~\eqref{eq:vinnicombe04controller} guarantees the same $\ell_2$-gain bound $\gamma^\star$ for both systems, however, this is now demonstrably suboptimal. He also showed that any controller that $\ell_2$-gain stabilizes both systems will achieve a closed-loop $\ell_2$-gain from $w$ to $x$ of at least $a/2$ for either of the systems. This has the following important implication: There exists, at least in discrete time, no $\ell_2$-gain stabilizing adaptive controller that is universal in the sense that it works for all $a\in\mathbb{R}$. Related work with suboptimal bounds on the $\ell_2$-gain from $w$ to $x$ has been published in~\cite{kjellqvist2022learning,kjellqvist2022minimax,orlov2018adaptive,rantzer2020ifac}.

\subsubsection{Lower bounds on the achievable $\ell_2$-gain}
Megretski~\cite{Megretski/R03} provided a lower-bound on the achievable $\ell_2$-gain from disturbances $w$ to states $x$ for scalar systems $x(t+1) = ax(t) + u(t) + w(t)$, where $a\in\mathcal{I}\subset\mathbb{R}$. Specifically, let $\gamma_0\geqslant 5$ and let $\mathcal{I}$ be a closed interval whose intersection with $\mathbb{R}_+$ or $\mathbb{R}_{-}$ has length at least $10\gamma_0$. Then any controller that achieves finite $\ell_2$-gain for all $a\in\mathcal{I}$, must necessarily have an $\ell_2$-gain from $w$ to $x$ greater than $\gamma_0$ for $a=0$. Thus, enlarging the model class $\mathcal{M}$ generally negatively impacts the $\ell_2$-gain from $w$ to $x$ even when the true system parameters are small.

\subsection{Standard zero-sum dynamic game formulation}
Next, we reformulate the adaptive optimal control problem introduced above as a standard zero-sum dynamic game, which can be addressed using dynamic programming.
\subsubsection{Compression of historic data}
Since the true system parameters $(A,B)\in\mathcal{M}$ are unknown, the optimal policy $\mu$ depends on past inputs $\{u(\tau)\}_{\tau=0}^{t-1}$ and states $\{x(\tau)\}_{\tau=0}^{t}$. That way, the optimal causal dual controller can ``learn'' the time-invariant system dynamics over time and adapt accordingly. Over time, keeping track of all these historic measurements becomes increasingly memory intensive and evaluation of the control policy grows more and more computationally expensive. To address this, the historical data is compressed in the matrix
\begin{equation}\label{eq:Zt}
    Z(t) = \sum_{\tau=0}^{t-1} \begin{bmatrix}
        x(\tau)\\
        u(\tau)\\
        x(\tau+1)
    \end{bmatrix}\begin{bmatrix}
        x(\tau)\\
        u(\tau)\\
        x(\tau+1)
    \end{bmatrix}^\top,
\end{equation}
whose dimensions remain constant over time. Similar (empirical) covariance matrices appear throughout the literature on learning for control, see, e.g.,~\cite{tsiamis2023statistical}. There is also a connection to the data matrices used in (off-line) data-driven control~\cite{Markovsky2008,Berberich2023,vanWaarde2020c,Meijer2025}, and subspace identification~\cite{vanOverschee1996} because $Z(t)$ is the product of a Hankel matrix and its transpose.

In addition to being computationally efficient, the compressed data $Z(t)$ satisfies, due to the time-invariance of $(A,B)\in\mathcal{M}$, the following important identity for all $t\in\mathbb{N}$:
\begin{equation}\label{eq:identity}
    \left\|\begin{bmatrix}
        A & B & -I
    \end{bmatrix}^\top\right\|^2_{Z(t)} = \sum_{\tau=0}^{t-1} \|w(\tau)\|^2,
\end{equation}
where $\|B\|_{A}^2=\trace B^\top AB$. As we will see next, Equation~\eqref{eq:identity} helps to reformulate the problem \eqref{eq:minimax-cost} as a standard zero-sum dynamic game~\cite{basar1999dynamic}, which can be addressed by dynamic programming. 

\begin{textbox}[!bt]
\sbsection{MINIMAX DYNAMIC PROGRAMMING}
\label{sb:minimax-dp}
The minimax optimal control problem can be reformulated as a standard zero-sum dynamic game, which can be addressed through minimax dynamic programming. Consider the minimax optimal control problem
\[
\begin{aligned}
    V_N(x_0)\coloneqq\inf_{\eta}\sup_{w}\quad
    & \sum_{t=0}^{N-1} g(x(t),u(t),w(t)) + g_N(x(N))\\
    \text{subject to}\quad
    & x(t+1) = f(x(t),u(t),w(t)),\quad u(t) = \mu_t(x(t)),
\end{aligned}
\]
where $x(t)\in\mathbb{R}^{n}$, $u(t)\in\mathbb{R}^{m}$, and $w(t)\in\mathbb{R}^{p}$ denote, respectively, the state, input, and disturbance at time $t\in\mathbb{N}$. The causal control policy $\mu=(\mu_0,\mu_1,\hdots,\mu_{N-1})$ generates inputs $u(t)$ based on past state measurements. The value function $V_N$ can be computed iteratively according to
\[
    V_{k+1}(x) = \min_u\max_w \left[g(x,u,w) + V_{k}(f(x,u,w))\right],\quad V_0(x) = g_N(x),
\]
where the minimizing $u$'s correspond to the optimal policy. If $\max_w g(x,u,w)\geqslant 0$ for all $x\in\mathbb{R}^{n}$ and $u\in\mathbb{R}^{m}$, then the sequence $\{V_k\}_{k=0}^{N}$ is monotonically increasing, i.e., $g_N(x) = V_{0}(x)\leqslant V_1(x)\leqslant\hdots\leqslant V_{N}(x)$ for all $x\in\mathbb{R}^n$. Moreover, $V_N$ is upper 
bounded by any function 
$\bar{V}$
satisfying
\[
    g_N(x) \leqslant \inf_u\sup_w\left[g(x,u,w) + \bar{V}(f(x,u,w))\right]\leqslant \bar{V}(x), 
\]
for all $x\in\mathbb{R}^{n}$. The latter provides the foundation for approximate dynamic programming and corresponds to the Bellman equation for infinite-horizon problems when equality holds. Such bounds can be computed by performing value iteration with the corresponding inequality. For details, we refer to, e.g.,~\cite{Wang2015,bertsekas2007dynamic,Bertsekas2012,basar1999dynamic}. 
\end{textbox}

\subsubsection{Reformulated problem}\label{sec:standard-game}
The textbox~``\nameref{sb:minimax-dp}'' summarizes basic dynamic programming for zero-sum dynamic games. These results are not directly applicable to the optimal dual control problem because the dynamics in Equation~\eqref{eq:system} depend on the unknown parameters $(A,B)\in\mathcal{M}$. However, Rantzer~\cite{rantzer2021minimax} reformulated the minimax optimal dual control problem so that standard dynamic programming can be applied. This reformulation is based on the fact that $w(t)$ enters the state dynamics additively, so the adversary can pick the next state $x(t+1)$ arbitrarily. Therefore, we define a new adversarial disturbance $v(t)\in\mathbb{R}^{n}$ so that $x(t+1) = v(t)$. Since $Z(t)$ is also the solution to a difference equation, the extended state $(x(t),Z(t))$ has the dynamics
\begin{equation}\label{eq:reformulated-dynamics}
    \begin{cases}
        x(t+1) &= v(t),\\
        Z(t+1) &= Z(t) + \begin{bmatrix}
        x(\tau)\\
        u(\tau)\\
        v(\tau)
    \end{bmatrix}\begin{bmatrix}
        x(\tau)\\
        u(\tau)\\
        v(\tau)
    \end{bmatrix}^\top, \quad Z(0)=Z_0=0.
    \end{cases}
\end{equation}
With Feldbaum's nomenclature, $Z(t)$ can be viewed as the information state, complementing the physical state $x(t)$.
Here, instead of putting $Z_0$ to zero, the initial condition can be exploited to incorporate prior knowledge about the system, e.g., in the form of previously collected data. 
Using the identity in Equation~\eqref{eq:identity}, we obtain the reformulated optimal control problem given by
\begin{equation}\label{eq:reformulated-cost}
    \hat{J}^\star(x_0,Z_0) \coloneqq \inf_{\eta}\sup_{v,N}\sup_{(A,B)\in\mathcal{M}}\sum_{t=0}^{N-1}\|x(t)\|^2_Q + \|u(t)\|^2_R - \gamma^2 \left\|\begin{bmatrix}
        A & B & -I
    \end{bmatrix}^\top\right\|_{Z(N)}^2,
\end{equation}
where $(x(t),Z(t))$ evolves according to Equation~\eqref{eq:reformulated-dynamics}, and the input is generated by a \emph{stationary} policy $u(t) = \eta(x(t),Z(t))$.
The resulting reformulated problem is a standard zero-sum dynamic game~\cite{basar1999dynamic}, which can be addressed using minimax dynamic programming. In fact, we will use dynamic programming in the next section to establish that the reformulated problem is equivalent to the original one. We will also see that, since $A$ and $B$ appear only in the terminal cost, the optimal policy $\eta$ depends only on the current state $(x(t),Z(t))$ and, in fact, this optimal policy is the stationary policy obtained through value iteration.

\subsection{Minimax dynamic programming}\label{sec:minimax-dynamic-programming}
Dynamic programming can be used to establish important equivalences between the original and the reformulated problem, as well as the corresponding optimal policies. In addition, it provides an approach to solving the reformulated problem, which, when combined with input randomization, has enabled the formulation of explicit dual control laws for certain important special cases of the adaptive minimax optimal control problem in Equation~\eqref{eq:minimax-cost}.

\subsubsection{Equivalence of the original and reformulated problems}
The reformulated problem above makes it natural to introduce the value iteration
\begin{equation}\label{eq:bellman-iteration}
    V_{k+1}(x,Z) = \min_u\max_v\left[\|x\|^2_Q + \|u\|^2_R + V_k\left(v,Z+(x,u,v)(x,u,v)^\top\right)\right],
\end{equation}
which is initialized by taking $V_0(x,Z)$ equal to the terminal cost in Equation~\eqref{eq:reformulated-cost}, i.e.,
\[
    V_0(x,Z) = -\gamma^2\min_{(A,B)\in\mathcal{M}}\left\|\begin{bmatrix}
        A & B & -I
    \end{bmatrix}^\top\right\|^2_Z.
\]
Standard dynamic programming arguments (see, e.g., the~``\nameref{sb:minimax-dp}'' textbox) are used to show that $V_N(x_0,0)$ is, for any $N\geqslant 0$, a lower bound for the value of the original problem in Equation~\eqref{eq:minimax-cost}, i.e., $V_N(x_0,0)\leqslant J^\star(x_0)$ for all $x_0\in\mathbb{R}^{n}$. Since the adversary can always choose $v=-(Ax+Bu)$ (or, equivalently, $w=0$ in the original formulation), the sequence $\{V_k(x,Z)\}_{k=0}^{\infty}$ grows monotonically. 
Equation~\eqref{eq:minimax-cost}~has a finite value if and only if the sequence $\{V_k(x,0)\}_{k=0}^{\infty}$ is upper bounded. If so, then the limit $V^\star\coloneqq \lim_{k\rightarrow\infty}V_k$ exists and is equal to the optimal value of the original problem, i.e.
\[
    V^\star(x_0,0) = J^\star(x_0)=\hat{J}^\star(x_0,0).
\]
Moreover, the optimal policy $\eta^\star$ for the reformulated problem is obtained as
\[
    \eta^\star(x,Z) = \arg\min_u\max_v\left[\|x\|^2_Q + \|u\|^2_R + V^\star(v,Z+(x,u,v)(x,u,v)^\top)\right].
\]
Note that this policy is stationary and depends only on the most recent value of the extended state $(x(t),Z(t))$. The optimal policy $\mu^\star$ for the original problem in Equation~\eqref{eq:minimax-cost} is given by
\[
    \mu^\star_t(x(0),\hdots,x(t),u(0),\hdots,u(t-1))\coloneqq \eta^\star(x(t),Z(t)).
\]
Since $V_k(x,Z)\leqslant V_k(x,0)$ for any $Z\succcurlyeq 0$, the limit then also exists for all $Z\succcurlyeq 0$ and prior information embedded in $Z$ may help to improve the performance. Finally,~\cite{rantzer2021minimax} also presented a suboptimal policy, obtained by approximating the value function, while extensions to state estimation and output feedback were derived in \cite{kjellqvist2024minimax}.


\begin{textbox}[!b]
{\small\sbsection{OPTIMAL EXPLORATION OVER THREE STEPS}\label{sb:optimal-exploration}
To provide some insight into the exploration--exploitation tradeoff, we will here go through the value iteration in Equation~\ref{eq:bellman-iteration}. We revisit the example from~\cite{vinnicombe2004examples} featuring the scalar system $x(t+1) = ax(t) + bu(t) + w(t)$, with sign uncertainty in $b$, i.e., $b\in\{-1,+1\}$, and known $a$. Let $Q=1$, $R=0$, and $\gamma \geqslant \gamma^\star$. The value iteration is initialized as
$
    V_0(x,Z) = -\gamma^2 \min_{b=\pm1}\left\|(a,b,-1)\right\|_{Z}^2.
$
We proceed according to Equation~\eqref{eq:bellman-iteration}: 
\begin{align}
    V_1(x,Z) &= \min_u\max_v \left[x^2 + V_0\left(v,Z+(x,u,v)(x,u,v)^\top\right)\right],\notag\\
    &= \min_u \left[x^2 - \gamma^2\min_v\min_{b=\pm1}\left(\left\|(a,b,-1)\right\|^2_Z +(ax+bu-v)^2\right)\right]
    = x^2 - \gamma^2\min_{b=\pm1}\left\|(a,b,-1)\right\|^2_Z.\notag\\\
    V_2(x,Z) &= \min_u\max_v\left[x^2 + V_1\left(v,Z+(x,u,v)(x,u,v)^\top\right)\right]\nonumber\\
    &= \max_\theta\min_u\max_v\sum_{b=\pm1}\theta_b\left[x^2 + v^2-\gamma^2(ax+bu-v)^2 - \gamma^2\left\|(a,b,-1)\right\|^2_Z \right],\label{eq:two-quadratics}\\
    &= \max_\theta\left(x^2 + \frac{a^2x^2}{1-\gamma^{-2}}\hbox{$\left[1-\left(\sum_{b=\pm1}b\theta_b\right)^2\right]$}
    -\gamma^2\sum_{b=\pm1}\theta_b\left\|(a,b,-1)\right\|^2_Z \right)\notag\\
    V_3(x,0) &= \min_u\max_{v,\theta}
    \left(x^2+v^2+\frac{a^2v^2}{1-\gamma^{-2}}\hbox{$\left[1-\left(\sum_{b=\pm1}b\theta_b\right)^2\right]$}
    -\gamma^2\sum_{b=\pm1}\theta_b(ax+bu-v)^2\right)\notag\\
    &= x^2+\min_u\max
    \left\{\frac{\gamma^2(\gamma^2-1+a^2\gamma^2)}{(\gamma^2-1)^2+a^2\gamma^2}a^2x^2
    -\gamma^2u^2,\frac{(ax+u)^2}{1-\gamma^{-2}},\frac{(ax-u)^2}{1-\gamma^{-2}}\right\}\label{eq:three-quadratics}
\end{align}
where the introduction of $\theta=(\theta_{-1},\theta_1)=(1-\theta_1,\theta_1)\in[0,1]^2$ makes it possible to shift the order of min and max. Notice that in Equation~\eqref{eq:two-quadratics}, we are minimizing the maximum of two convex quadratic functions, while in Equation~\eqref{eq:three-quadratics}, we are minimizing three quadratic functions, where one of them is concave. This is illustrated in Figure~\ref{fig:minimax-illustration}. The concave function activates the input leading to exploration. }
\end{textbox}
\begin{figure}[!ht]
    \centering
        
            





            
    \begin{minipage}{.45\textwidth}
        \centering
        \begin{tikzpicture}[scale=2, transform shape]
            \pgfmathsetmacro{\a}{1} 
            \pgfmathsetmacro{\gamma}{3} 
            \pgfmathsetmacro{\xi}{.5} 
            \pgfmathsetmacro{\ginv}{1/(1-(\gamma)^(-2)} 
            \pgfmathsetmacro{\offset}{(1+\ginv*(\a)^2)*(\xi)^2} 
            \pgfmathsetmacro{\tone}{1+\ginv*(\a)^2} 
            \pgfmathsetmacro{\ttwo}{1+((\a)^2)/((1/\tone)-(\gamma)^(-2))}
        
            \useasboundingbox (-1,-.3) rectangle (1,2);
            \clip (-1,-.3) rectangle (1,2);
            
            \draw[->] (-1,0) -- (1,0) node[below,pos=0.95,font=\tiny] {$u$};
            \draw[->] (0,-.3) -- (0,2) node[left] {};

            \draw[name path=A,thick,blue]
            plot[domain=-1:1,samples=300]
            (\x,{\offset + \ginv*(\x*(\x-2*\a*\xi))});
            \draw[name path=B,thick,blue]
            plot[domain=-1:1,samples=300]
            (\x,{\offset + \ginv*(\x*(\x+2*\a*\xi))});

            \path[name intersections={of=A and B, by=I1}];


            \draw[name path=C,thick,red]
            plot[domain=-1:1,samples=300]
            (\x,{\ttwo*(\xi)^2 - (\gamma)^2*(\x)^2});

            \path[name intersections={of=A and C, by=I2}];
            \path[name intersections={of=B and C, by={Q,I3}}];
            
            \fill[red] (I2) circle (1.5pt);
            \fill[red] (I3) circle (1.5pt);
        \end{tikzpicture}\\[2mm]
        {\small {\bf(a)} Optimal actuation minimizes the maximum of three parabolas.}
    \end{minipage}\hfill    
    \begin{minipage}{.45\textwidth}
        \centering
        \begin{tikzpicture}[scale=2, transform shape]
            \pgfmathsetmacro{\a}{1} 
            \pgfmathsetmacro{\gamma}{3} 
            \pgfmathsetmacro{\xi}{.5} 
            \pgfmathsetmacro{\ginv}{1/(1-(\gamma)^(-2)} 
            \pgfmathsetmacro{\offset}{(1+\ginv*(\a)^2)*(\xi)^2} 
            \pgfmathsetmacro{\tone}{1+\ginv*(\a)^2} 
            \pgfmathsetmacro{\ttwo}{1+((\a)^2)/((1/\tone)-(\gamma)^(-2))}
        
            \useasboundingbox (-1,-.3) rectangle (1,2);
            \clip (-1,-.3) rectangle (1,2);
            
            \draw[->] (-1,0) -- (1,0) node[below,pos=0.95,font=\tiny] {$u$};
            \draw[->] (0,-.5) -- (0,2) node[left] {};

            \draw[name path=MAX,thick,blue]
            plot[domain=-1:1,samples=300]
            (\x,{max(\offset + \ginv*(\x*(\x-2*\a*\xi)),\offset + \ginv*(\x*(\x+2*\a*\xi)),\ttwo*(\xi)^2 - (\gamma)^2*(\x)^2)});

            \pgfmathsetmacro{\qA}{\ginv+(\gamma)^2}
            \pgfmathsetmacro{\qB}{2*\ginv*\a*\xi}
            \pgfmathsetmacro{\qC}{\offset-\ttwo*(\xi)^2}
            \pgfmathsetmacro{\qDisc}{(\qB)^2-4*\qA*\qC}
            \pgfmathsetmacro{\xstar}{(-\qB+sqrt(\qDisc))/(2*\qA)}
            \pgfmathsetmacro{\ystar}{\ttwo*(\xi)^2-(\gamma)^2*(\xstar)^2}
            
            \coordinate (Iright) at (\xstar,\ystar);
            \coordinate (Ileft) at (-\xstar,\ystar);
            
            \fill[red] (Ileft) circle (1.5pt);
            \fill[red] (Iright) circle (1.5pt);
        \end{tikzpicture}\\[2mm]
        {\small {\bf(b)} Non-convex objective}
    \end{minipage}\\[2mm]
    \caption{After three steps, value iteration requires minimizing the maximum of two convex (\kern0.5pt\protect\tikz[scale=0.6]{\protect\draw[color=blue,line width=1pt] (0,0.1)--(17pt,0.1);\protect\draw[color=white] (0,0)--(0.55,0)}) and one concave (\kern0.5pt\protect\tikz[scale=0.6]{\protect\draw[color=red,line width=1pt] (0,0.1)--(17pt,0.1);\protect\draw[color=white] (0,0)--(0.55,0)}) parabola to determine the optimal actuation. As can be seen in {\bf (a)}, the concave parabola shifts the maximizing $u$ away from zero, causing the controller to probe the system by applying non-zero actuation even when no prior knowledge is available ($Z=0$). However, it also renders the objective (\kern0.5pt\protect\tikz[scale=0.6]{\protect\draw[color=blue,line width=1pt] (0,0.1)--(17pt,0.1);\protect\draw[color=white] (0,0)--(0.55,0)}) non-convex, as shown in {\bf (b)}. Input randomization helps to remedy this non-convexity.}
    \label{fig:minimax-illustration}
\end{figure}

The textbox~``\nameref{sb:optimal-exploration}'' illustrates how the value iteration in Equation~\ref{eq:bellman-iteration} over three steps naturally leads to an exploration--exploitation tradeoff. It also shows that the minimization problem in the Bellman recursion eventually becomes non-convex, which makes it difficult to solve. However, as we will see in the next subsection, introduction of randomized control policies helps to overcome this challenge and recover the convexity.

\subsubsection{Randomized control policies}
The value iteration simplifies dramatically if we can switch the order of minimization and maximization. By Von Neumann's minimax theorem~\cite{neumann1928theorie}, this is possible for convex-concave objectives. Specifically, we use the following formulation by Fan~\cite{fan1953minimax}, which was communicated by von Neumann in 1952: Let $X\subset\mathbb{R}^{n}$ and $Y\in\mathbb{R}^{m}$ be compact convex sets. If $f:X\times Y\rightarrow \realR$ is a concave-convex continuous function,
then
\begin{eqnarray*}
    \max_{x\in X}\min_{y\in Y}f(x,y)=\min_{y\in Y}\max_{x\in X}f(x,y).
\end{eqnarray*}
Applying input randomization helps to obtain convexity in $u$. This means that at time $t$ the policy $\mu$ generates a random input vector $u(t)$ (with finite variance). The cost function is modified to optimize the expected value, i.e.,
\begin{equation*}
    J^\star(x_0)\coloneqq \inf_{\mu}\sup_{w,N}\sup_{(A,B)\in\mathcal{M}}
    \averageE\sum_{t=0}^{N-1}\left(\|x(t)\|^2_Q + \|u(t)\|^2_R - \gamma^2\|w(t)\|^2\right).
\end{equation*}
To see how this helps, consider $\mathcal{M}=\{(A_\theta,B_\theta):\theta\in\Theta\}$ with the symmetry property that $\Theta=-\Theta$ and $(A_{-\theta},B_{-\theta})=(A_{\theta},-B_{\theta})$. Every model instance $(A_\theta,B_\theta)$ is supposed to be stabilizable with a corresponding Lyapunov function $|x|^2_{P_\theta}$ satisfying
\[
    |x|^2_{P_\theta}\geqslant\min_u\max_v\left(|x|^2_Q+|u|^2_R-\gamma^2|A_\theta x+B_\theta u-v|^2+|v|^2_{P_\theta}\right).
\]
In many cases, the argument in the previous sections can be generalized to show that the Bellman inequality
\begin{align}
	\bar{V}(x,Z)\geqslant\min_u\max_v\left[\averageE\left(|x|^2_Q+|u|^2_R\right)+\bar{V}\left(v,Z+(x,u,v)(x,u,v)^\top\right)\right]
	\label{eqn:Bellmanineq}
\end{align}
has explicit solutions of the form $\bar{V}(x,Z)\coloneqq \max_{\theta}\max\{V_{\theta}(x,Z),W_{\theta}(x,Z)\}$, with
\begingroup
\begin{align*}
    V_{\theta}(x,Z)&= \averageE|x|^2_{P_\theta}
    -\phantom{\Bigg|}\gamma^2\left\|
    \begin{bmatrix}A_\theta&B_\theta&-I\end{bmatrix}^\top\right\|^2_{\mathbb{E}Z},\\
    W_{\theta}(x,Z) &= \averageE|x|^2_{T_\theta}-\frac{\gamma^2}{2}\sum_{\mathclap{k=\pm1}}\left\|
    \begin{bmatrix}A_\theta&kB_\theta&-I\end{bmatrix}^\top\right\|^2_{\mathbb{E}Z}.
\end{align*}
\endgroup
Here, $|x|^2_{P_\theta}$ quantifies the expected cost with a known model $(A_\theta,B_\theta)$, and $|x|^2_{T_\theta}$ is the expected cost for the same model when it is not known. 
Note that the bracket on the right hand side of the inequality \eqref{eqn:Bellmanineq} is a maximum of expressions that are linear in $(\averageE u,\averageE u^2,Z)$. After maximization over $v$ and minimization over $u$, the right hand side is still convex in $Z$. To verify the inequality \eqref{eqn:Bellmanineq}, it is therefore sufficient to focus on the extreme values of $Z$, with $Z=0$ (no data) being critical. (Without randomized inputs, this would not have been the case.)
Once the inequality holds, a stabilizing dual control law mapping $(x,Z)$ to a random $u$ is obtained by solving the convex optimization problem over $(\averageE u,\averageE u^2)$ on the right hand side.


Details of this argument have been worked out in~\cite{rantzer2026minimax} for single-input systems with a known $A$-matrix. The $B$-vector is unknown, but stabilizability and boundedness are assumed in the sense that
\begin{align*}
    \min_u\max_w\left(|x|^2_Q+|u|^2_R-\gamma^2|w|^2+|Ax+Bu+w|^2\right)\leqslant|x|^2\text{ and }
			1\leqslant|B|\leqslant \beta.
\end{align*}
The resulting minimax dual controller computes, at each time instant $t$, a worst-case value of $\hat{B}_t$ based on $(x(t),Z(t))$ using a semi-definite program. Let $K_t$ be the minimax optimal state-feedback gain for the model $(A,\hat{B}_t)$. Then, the optimal controller has the form
\begin{equation}\label{eq:single-input-minimax}
\begin{cases}
    u(t) = K_tx(t)\qquad &\text{if past data is sufficiently rich}\\
    \mathbb{E}u(t)=\sigma K_tx(t),~\mathbb{E}[u^2(t)] = |K_tx(t)|^2,\qquad &\text{otherwise},
\end{cases}
\end{equation}
where the number $\sigma\in[0,1]$ quantifies the confidence in the value $\hat{B}_t$ as opposed to $-\hat{B}_t$. As expected from the construction, the control action depends continuously on the state: Without data, $\sigma=0$ and, thus, $\mathbb{E}u(t)=0$, but as more data becomes available $\sigma$ grows and eventually as $\sigma=1$ the controller becomes a certainty-equivalence controller. 



\section{NUMERICAL CASE STUDIES}
In this section, we present numerical case studies in which we apply some of the different adaptive/dual control approaches discussed before. We consider the ARMAX model
\begin{equation}\label{eq:ARMAX}
    y(t+1) = -a_0y(t) - a_1y(t-1) + b_0u(t) + b_1u(t-1) + e(t+1) + c_0e(t) + c_1e(t-1).
\end{equation}
where the coefficients $a_0$ and $a_1$ are \emph{known}, while $b_0$, $b_1$, $c_0$, and $c_1$ are \emph{unknown}. Here, $e(t)$ is Gaussian white noise. This system admits a state-space realization
\[
    x(t+1) = \begin{bmatrix}
        -a_0 & 1\\
        -a_1 & 0
    \end{bmatrix}x(t) + \begin{bmatrix}
        b_0\\
        b_1
    \end{bmatrix}u(t) + w(t),\qquad y(t) = Cx(t),
\]
where the noise is absorbed into $w(t)$ by taking $w(t) = (c_0e(t)+e(t+1),c_1e(t))$. However, note that the additive disturbance $w(t)$ is more general and can also be used to inject, e.g., $n$-dimensional Gaussian white noise or worst-case disturbances to investigate the settings of~\cite{jedra2022minimal} and~\cite{rantzer2026minimax} respectively. This will be done in the sequel to compare the performance of several different adaptive/dual control schemes. To meet the different assumptions of these schemes, it may be needed to transform the above state-space realization to an appropriate coordinate system. While Equation~\eqref{eq:ARMAX} would lead us to take $C=\begin{bmatrix} 1 & 0\end{bmatrix}$, we assume full state information ($C=I$) because the state consists of past inputs and outputs and doing so allows the methods~\cite{rantzer2026minimax,jedra2022minimal} to be investigated. In the above realization, the $A$-matrix is known while $B$ is unknown. Uncertainty in the $B$-matrix plays a central role in dual control because, unlike the $A$-matrix, $B$ cannot be studied without activation of the input.


\subsection{FIR model with Gaussian white noise}
First, we take $a_0 = a_1 = 0$ and let the true parameters be $b_0 = 1$ and $b_1 = 1.5$. We compare the following controllers: The minimax controller in Equation~\eqref{eq:single-input-minimax} by Rantzer~\cite{rantzer2026minimax}, the regret rate minimizing certainty-equivalence controller (CEC) by Jedra \& Proutiere~\cite{jedra2022minimal}, and the certainty-equivalence linear-quadratic regulator based on weighted recursive least-squares (wRLS CE-LQR) by Guo~\cite{guo1996self}. It should be noted that the system is non-minimum phase, i.e., the $B$-polynomial has roots outside the unit circle, which is traditionally challenging to deal with for minimum-variance self-tuning regulators as discussed in Section~\ref{sec:STR}. All of the controllers considered here, however, are able to deal with such non-minimum phase behavior. For each of the controllers, we use a forgetting factor of $0.98$, $Q=0.1I$, $R=0.01$, and, for the minimax approach, $\gamma=10$ and $\|B\|\leqslant \beta=3$. 

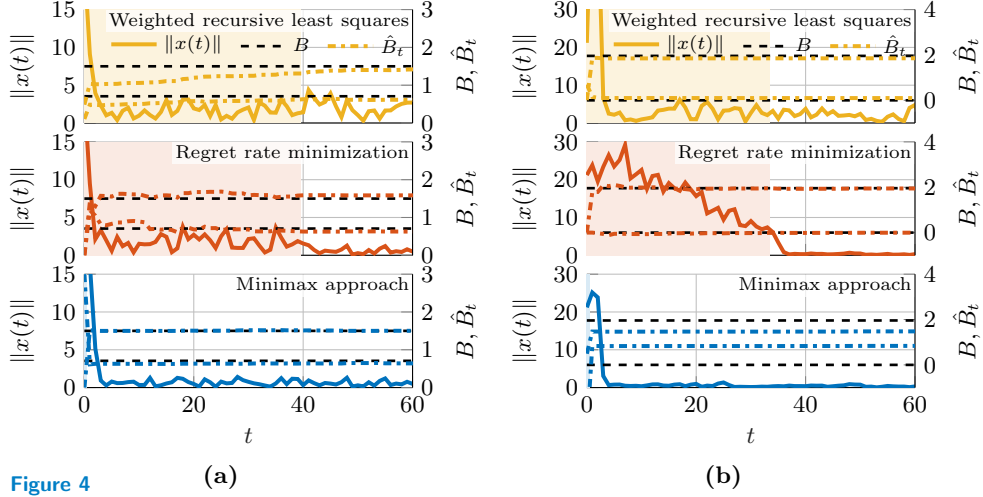
\begin{figure}[!bt]
\centering
    \hspace*{-.75cm}
    \begin{minipage}{.45\linewidth}
        \centering
        \input{figures/nmp-gaussian_auto}
        {\small {\bf(a)}}
    \end{minipage}\hspace*{.75cm}   
    \begin{minipage}{.45\linewidth}
        \centering
        \input{figures/Ljung-colored_auto}
        {\small {\bf(b)}}
    \end{minipage}
    \caption{Simulation results for {\bf (a)} $a_0 =a_1=0$, $b_0=1$, $b_1=1.5$ and zero-mean Gaussian white noise $w$, and {\bf (b)} $a_0 = 0.9$, $a_1=0.95$, $b_0=10$, $b_1=0$ and colored noise $w$ ($c_0 = 2c_1= 1.5$). Top: Weighted recursive least squares~\cite{guo1996self}. Middle: Regret rate minimization~\cite{jedra2022minimal}. Bottom: Minimax approach~\cite{rantzer2026minimax}. Each panel uses two y-axes with $\|x(t)\|$ on the left, and the true parameters $B$ and their estimates $\hat{B}_t$ on the right. The shaded regions indicate when the controller is actively probing the system. \textbf{Caveat:} There is no doubt that more extensive analysis and simulations would be needed before drawing any general conclusions about applicability of the different approaches.}
    \label{fig:simulation-results_v2}
\end{figure}
Figure~\ref{fig:simulation-results_v2}(a) shows the resulting trajectories when the disturbance $w$ is zero-mean Gaussian white noise with variance $\sigma^2=0.25$. It can be seen that all three controllers are stabilizing despite the non-minimum-phase behavior. Interestingly, it turns out that the minimax controller only performs exploration in the first time-step. This is due to the fact that in the scalar input case, the parameter vector $B$ has the same dimension as the state. Thus, a single exploration step is sufficient to receive information about all $B$-parameters.
Both the CEC and the wRLS CE-LQR explore over a longer time span. Consequently, their parameter estimates $\hat{B}_t$ converge much faster, at the cost of (short-term) control performance.

\subsection{ARMAX model with colored noise} We compare again the same controllers but now, inspired by~\cite[Example 1]{ljung1975counterexamples}, we take $a_0 = 0.9$ and $a_1=0.95$, and $b_0=10$ and $b_1 = 0$ for the true parameters. In contrast with the previous example, we now inject colored noise by setting $c_0 = 2c_1 = 1.5$. While the minimax controller should obviously be able to handle colored noise, this represents a significant deviation from the Gaussian white noise setting considered in~\cite{jedra2022minimal}, and the condition that Equation~\eqref{eq:Cinv_min_half} is positive real does not hold contrary to the assumptions in~\cite{guo1996self}. Interestingly, the estimates for CEC and wRLS CE-LQR still converge as visible in Figure~\ref{fig:simulation-results_v2}(b), where we also see that, although all controllers are stabilizing, the minimax approach drives the state close to zero much faster than the other controllers. Interestingly, this is despite the fact that the estimates $\hat{B}_t$ of the minimax approach do not converge to the exact values. Nevertheless, the performance is not significantly affected by this and, therefore, the minimax controller still only performs a single exploration. The CEC and wRLS CE-LQR schemes again perform exploration over a longer time span, which, especially for CEC, sacrifices exploitation. After the exploration phase ends, however, CEC rapidly reaches comparable performance to the minimax approach. 

\section{CONCLUSIONS}
The exploration--exploitation tradeoff has been studied extensively for almost a century and plays a crucial role across many scientific disciplines. Learning-based control is no different: Optimal controllers must balance immediate performance against the value of new information. 

The most complete theory for exploration--exploitation tradeoffs has been attained for the multi-armed bandit problem. This includes important concepts such as the Gittins index, ``optimism in the face of uncertainty'', and regret minimization. In the field of adaptive control, self-tuning regulators have demonstrated that remarkably simple certainty-equivalence ideas can be highly effective in practice, and they have also revealed important theoretical issues related to, e.g., persistence of excitation and closed-loop identifiability. More recently, regret analysis for linear-quadratic control has clarified when exploration is fundamentally necessary and has led to learning algorithms with provably optimal (asymptotic) regret rate. Finally, recent developments of the minimax dual control paradigm have resulted in policies with strong non-asymptotic performance guarantees and provable robustness to unmodelled dynamics. 

Despite this truly substantial progress, many open challenges and opportunities remain. In particular, the minimax approach requires further research to deal with more general and challenging settings. Another important direction for future research is to handle high-dimensional and data-scarce systems, where scalability of the control algorithm becomes an important issue and where too few measurements to estimate all relevant system parameters. Again, a properly-made exploration--exploitation tradeoff should still allow the most relevant system parameters to be estimated so that satisfactory performance can be attained. A promising approach for addressing scalability is to exploit structure of positive systems, see, e.g.,~\cite{RantzerValcher2021}. Recently, the Bellman equation was solved explicitly for positive systems with linear cost~\cite{rantzer2022explicit}, which significantly broadens the types of (dual) control problems for which we hope to obtain exact optimal solutions. In fact, it has already been extended to the minimax setting~\cite{gurpegui2023minimax} and found applications to adaptive control~\cite{bencherki2025}, making minimax dual control for positive systems a natural next step. A final direction that deserves further research is the effect of including non-zero prior information in the initial condition $Z_0$ in the minimax approach. The minimax dual control problem naturally has a finite-horizon counterpart, where this initial condition can be used to solve the problem in a receding-horizon fashion. This provides a compelling path toward integration of the minimax dual control approach with model predictive control. 

\begin{summary}[SUMMARY POINTS]
\begin{enumerate}
\item 
Bayesian formulations of {\bf multi-armed bandit problems were completely solved} in terms the Gittins index. A frequentist formulation was solved for optimal asymptotic regret rate using a policy based on ``optimism in the face of uncertainty''.

\item 
Their algorithmic simplicity, flexible disturbance models, and {\bf extensively-documented industrial applications make self-tuning controllers invaluable benchmarks} for any dual control algorithm.

\item Over the past twenty years, regret has been the dominant quality measure for online learning in linear-quadratic control, and {\bf it has led to deeper understanding of what makes a dual control problem hard.} For problems where self-tuning control belongs, logarithmic growth of regret is achievable. Other problems require extensive exploration near the optimal policy, resulting in a faster growth of regret.

\item In the {\bf minimax dual control} paradigm, randomized control policies allow explicit solutions to the Bellman inequality to be obtained, leading to policies with truly non-asymptotic bounds on the cost and provable robustness to disturbances that are correlated over time or generated by unmodelled non-linear dynamics.
\end{enumerate}%
\end{summary}

\begin{issues}[FUTURE ISSUES]
\begin{enumerate}
\item The minimax paradigm has led to promising new dual controllers, but further research is required to extend these results to general uncertain linear systems. Furthermore, regret analysis for the minimax approach would provide valuable insight into how its performance compares to other dual control methods.
\item Exploit additional structural properties, such as, e.g., positivity and sparsity, of the system to develop scalable and computationally-efficient dual control algorithms.
\item In high-dimensional systems, relatively few measurements are available compared to the number of unknown system parameters. The development and analysis of dual controllers that deal with such high-dimensional but data-scarce systems is challenging but important for many modern control applications.
\item Integrate dual control principles with constrained control techniques to deal with actuator or safety/performance constraints.
\end{enumerate}%
\end{issues}

\section*{DISCLOSURE STATEMENT}
The authors are not aware of any affiliations, memberships, funding, or financial holdings that
might be perceived as affecting the objectivity of this review. 

\section*{ACKNOWLEDGMENTS}
Support was received from the European Research Council (Advanced Grant 101199738) and Wallenberg AI, Autonomous Systems and Software Program (WASP). Views and opinions expressed are however those of the author(s) only and do not necessarily reflect those of the European Union or the European Research Council (ERC) Executive Agency. Neither the European Union nor the granting authority can be held responsible for them. 

The authors are members of the Excellence Center at Link\"{o}ping-Lund in Information Technology (ELLIIT). 

\bibliographystyle{ar-style3}
\bibliography{references}

\end{document}

%% file: figures/nmp-gaussian_auto.tex
\definecolor{mycolor1}{rgb}{0.00000,0.44700,0.74100}
\definecolor{mycolor2}{rgb}{0.85000,0.32500,0.09800}
\definecolor{mycolor3}{rgb}{0.92900,0.69400,0.12500}
\definecolor{mycolor4}{rgb}{0.12941,0.12941,0.12941}%

\newlength{\tikzplotheight}
\newlength{\vplotspacing}
\newlength{\tikzplotwidth}
\setlength{\tikzplotheight}{1.5cm}
\setlength{\vplotspacing}{.25cm}
\setlength{\tikzplotwidth}{.75\linewidth}

\def\exploreEndCEC{39}
\def\exploreEndWRLS{39}
\def\exploreEndMinimax{0}

\def\xMinAll{0}
\def\xMaxAll{60}
\def\yMinNorm{0}
\def\yMaxNorm{15}
\def\yMinB{0}
\def\yMaxB{3}

\begin{tikzpicture}

\pgfplotsset{
  vaxisbase/.style={
    width=\tikzplotwidth,
    height=\tikzplotheight,
    scale only axis,
    xmin=\xMinAll, xmax=\xMaxAll,
    xtick={0,20,40,60,80,100,120,140},
    axis background/.style={fill=none},
    clip=true,
  },
}

\begin{axis}[
  vaxisbase,
  at={(0,2\tikzplotheight+2\vplotspacing)},
  axis x line*=bottom,
  axis y line*=left,
  xticklabels={},
  ymin=\yMinNorm, ymax=\yMaxNorm,
  ylabel={$\|x(t)\|$},
  legend style={legend cell align=left, align=left, font=\scriptsize,
    at={(0.05,0.87)}, anchor=north west, draw=none, fill=none, /tikz/every even column/.append style={column sep=4pt}},
  ymajorgrids, xmajorgrids,
]
\fill[mycolor3!15,draw=none]
  (axis cs:-0.5,\yMinNorm) rectangle (axis cs:\exploreEndWRLS+0.5,\yMaxNorm);
\addplot [color=mycolor3, line width=1.5pt] table[row sep=crcr]{%
0	21.2132034355964\\
1	9.11423681684283\\
2	3.95856683403309\\
3	2.32510966840068\\
4	1.01882180674386\\
5	1.34309612641006\\
6	0.385305353770782\\
7	1.61344782340527\\
8	2.01822097342086\\
9	0.598264210265117\\
10	0.967493415374335\\
11	1.66113168734257\\
12	1.05101847284035\\
13	0.594850661463536\\
14	1.4533953693466\\
15	3.04329492902101\\
16	1.06014908206893\\
17	1.60348690333643\\
18	3.09300473760634\\
19	3.11852304366506\\
20	1.92754357280447\\
21	1.71639440056148\\
22	2.38359986406556\\
23	0.908195367053798\\
24	2.38461863631837\\
25	0.467298941210309\\
26	0.945314914237159\\
27	0.570632532833491\\
28	1.40148482994172\\
29	1.76519041178866\\
30	2.05258515882702\\
31	0.400927198846367\\
32	2.70707506663944\\
33	2.28767884883605\\
34	1.38585134498777\\
35	1.36010627259026\\
36	3.0740543422079\\
37	2.21872619835528\\
38	1.03865979355738\\
39	1.28640716720743\\
40	1.41175763622161\\
41	4.3574858341967\\
42	3.30835938681773\\
43	2.44872397122163\\
44	2.77760215217825\\
45	3.91071000678477\\
46	2.21077435962994\\
47	0.770062203154915\\
48	3.00911521885343\\
49	1.70269011892572\\
50	1.042888962747\\
51	0.294806840171242\\
52	1.55654576360823\\
53	1.6528143108858\\
54	0.570816724634402\\
55	1.65067786881898\\
56	1.66432634633532\\
57	2.18742678796663\\
58	2.5289058150552\\
59	2.68768992056339\\
60	2.70690240571913\\
61	2.24209819258412\\
62	1.73612586652971\\
63	0.965282679994246\\
64	1.14000317551535\\
65	1.6441420709727\\
66	0.145633266292823\\
67	1.1981060057787\\
68	1.08416323600987\\
69	1.3006467263663\\
70	0.673665076957973\\
71	2.38036611844986\\
72	2.57972368151228\\
73	0.405788879726043\\
74	2.08862065230708\\
75	1.39620526994057\\
76	1.72857420213411\\
77	1.08231831065021\\
78	1.92683754111655\\
79	1.00300765845943\\
80	1.28086496255089\\
81	1.48357594789464\\
82	1.4477258701032\\
83	1.57805109180883\\
84	2.11173924812789\\
85	2.90660748245581\\
86	1.4659549240794\\
87	1.36390593947296\\
88	0.316193335448858\\
89	0.72875301891704\\
90	0.897761726427897\\
91	0.466673774572479\\
92	0.554884585703448\\
93	0.415015862614084\\
94	1.74893333129195\\
95	1.94837509494526\\
96	1.10969266761898\\
97	0.443689071811988\\
98	2.02901788997482\\
99	2.33853664428133\\
100	2.03406542092715\\
};
\addlegendentry{$\|x(t)\|$}
\node[anchor=north east,font=\scriptsize,fill=white,fill opacity=0.75,text opacity=1,inner sep=1pt] at (axis cs:\xMaxAll,\yMaxNorm) {Weighted recursive least squares};
\end{axis}

\begin{axis}[
  vaxisbase,
  at={(0,2\tikzplotheight+2\vplotspacing)},
  axis x line=none,
  axis y line*=right,
  ymin=\yMinB, ymax=\yMaxB,
  ylabel={$B,\hat B_t$},
  ylabel style={color=black},
  legend style={legend cell align=left, align=left, font=\scriptsize,
    at={(1.02,0.87)}, anchor=north east, draw=none, fill=none, /tikz/every even column/.append style={column sep=4pt}},
  legend columns=2,
]
\addplot [color=black, dashed, line width=1pt, forget plot] table[row sep=crcr]{%
0	0.707106781186547\\
1	0.707106781186547\\
99	0.707106781186547\\
};
\addplot [color=black, dashed, line width=1pt] table[row sep=crcr]{%
0	1.5\\
1	1.5\\
99	1.5\\
};
\addlegendentry{$B$}
\addplot [color=mycolor3, line width=1.5pt, dashdotted] table[row sep=crcr]{%
0	0.106025470433517\\
1	0.483834653286145\\
2	0.484601929208136\\
3	0.469530827455097\\
4	0.475957844412999\\
5	0.482422668959872\\
6	0.483185624439942\\
7	0.491917242558874\\
8	0.49991435103167\\
9	0.499138906501349\\
10	0.499510515078539\\
11	0.507047620506253\\
12	0.522971935484304\\
13	0.523898019329849\\
14	0.528263611449646\\
15	0.545726061029557\\
16	0.541387227788664\\
17	0.542579000651582\\
18	0.570899647633059\\
19	0.577120076887667\\
20	0.577843762220822\\
21	0.58344603243684\\
22	0.559310580587519\\
23	0.567372537227684\\
24	0.572276867184506\\
25	0.580146570592807\\
26	0.5839735830893\\
27	0.583218348187911\\
28	0.583101850852454\\
29	0.584870750994201\\
30	0.587673658055519\\
31	0.587539595450631\\
32	0.597954040354293\\
33	0.596941332226205\\
34	0.596411951681517\\
35	0.595058520945827\\
36	0.585315375215314\\
37	0.606427512680326\\
38	0.618063600560395\\
39	0.61788622420229\\
40	0.628826344601178\\
41	0.62693091461462\\
42	0.626378954534319\\
43	0.635015458491024\\
44	0.638486528063179\\
45	0.636360517179854\\
46	0.636226241379775\\
47	0.635884169063655\\
48	0.611238431324741\\
49	0.615052919218807\\
50	0.615206309529796\\
51	0.615203821487817\\
52	0.609469831062895\\
53	0.610554637081398\\
54	0.610401270872843\\
55	0.614935934698002\\
56	0.615516284388639\\
57	0.614189852280307\\
58	0.616374974328424\\
59	0.60754856250271\\
60	0.603085213147391\\
61	0.605169918500922\\
62	0.607796753178932\\
63	0.607422259724077\\
64	0.608455480993925\\
65	0.607976983473942\\
66	0.608461094182346\\
67	0.605407416943149\\
68	0.60515038636255\\
69	0.604728103951462\\
70	0.60137306352992\\
71	0.608882796211212\\
72	0.613488607943922\\
73	0.614876637753463\\
74	0.615925345383572\\
75	0.61736928181288\\
76	0.618989825522722\\
77	0.62159164230415\\
78	0.619942748703425\\
79	0.621735363858007\\
80	0.623051561009972\\
81	0.619696398601166\\
82	0.619512065456643\\
83	0.619294448459274\\
84	0.611739057753909\\
85	0.615597572286616\\
86	0.615488317704842\\
87	0.615175718096471\\
88	0.613190144506007\\
89	0.612792691939024\\
90	0.612472490985632\\
91	0.610045409193183\\
92	0.609678578124538\\
93	0.610427709324356\\
94	0.612572719803695\\
95	0.613960024424902\\
96	0.612806587226967\\
97	0.613265319130762\\
98	0.613158331381422\\
99	0.614763976966271\\
};
\addlegendentry{$\hat{B}_t$}
\addplot [color=mycolor3, line width=1.5pt, dashdotted, forget plot] table[row sep=crcr]{%
0	0.484799324729201\\
1	1.0358020414653\\
2	1.02587430534113\\
3	1.02484700695948\\
4	1.0324720190116\\
5	1.0338490910193\\
6	1.03527846800791\\
7	1.05245169390605\\
8	1.05731048977289\\
9	1.0584708356794\\
10	1.05996323123244\\
11	1.06450426073453\\
12	1.07361220922723\\
13	1.07490024782688\\
14	1.08938048495472\\
15	1.10902246192684\\
16	1.10909386707727\\
17	1.10933379469385\\
18	1.16209365136322\\
19	1.17172286481139\\
20	1.19737611074112\\
21	1.20678859811773\\
22	1.22631076426306\\
23	1.22805040619567\\
24	1.23849831203398\\
25	1.23829270051751\\
26	1.2398043116582\\
27	1.23776653045495\\
28	1.24231962304534\\
29	1.2419975569979\\
30	1.24904452586098\\
31	1.24912675685008\\
32	1.27381385663317\\
33	1.27580116346537\\
34	1.27520903981246\\
35	1.27870875741742\\
36	1.30975943092445\\
37	1.30782769770953\\
38	1.30821928998024\\
39	1.3094698990382\\
40	1.30933486636655\\
41	1.3488323553348\\
42	1.35096786160571\\
43	1.36523051678467\\
44	1.36681716543137\\
45	1.38602336707973\\
46	1.38518998302041\\
47	1.385148990078\\
48	1.3993723469703\\
49	1.40161744477528\\
50	1.40184454072207\\
51	1.40184342971793\\
52	1.40160339774897\\
53	1.40184225538117\\
54	1.4002265146207\\
55	1.40019345410151\\
56	1.39775987260766\\
57	1.40286943622073\\
58	1.40528640788992\\
59	1.41373218584953\\
60	1.41836387370491\\
61	1.4243208923078\\
62	1.42828102190875\\
63	1.42857341384418\\
64	1.42908465617698\\
65	1.43323967980942\\
66	1.43291717504387\\
67	1.43388967557073\\
68	1.43353407312258\\
69	1.43479599103026\\
70	1.43205018503855\\
71	1.43302688401262\\
72	1.43344682822739\\
73	1.43309211963409\\
74	1.43521483177372\\
75	1.43718991994624\\
76	1.43773668198156\\
77	1.43953995736196\\
78	1.44204653284789\\
79	1.44337263014428\\
80	1.43822914723913\\
81	1.44128520197607\\
82	1.44117494488706\\
83	1.44534905326527\\
84	1.44547131572653\\
85	1.4494143218265\\
86	1.44965774536489\\
87	1.44898846924146\\
88	1.44817992819602\\
89	1.44807156405787\\
90	1.44917539578671\\
91	1.44882147624163\\
92	1.44866917952627\\
93	1.44778790862655\\
94	1.45169635280671\\
95	1.45096604787761\\
96	1.45214349346878\\
97	1.45205740506032\\
98	1.45586243919171\\
99	1.46229801320241\\
};
\end{axis}

\begin{axis}[
  vaxisbase,
  at={(0,\tikzplotheight+\vplotspacing)},
  axis x line*=bottom,
  axis y line*=left,
  xticklabels={},
  ymin=\yMinNorm, ymax=\yMaxNorm,
  ylabel={$\|x(t)\|$},
  ymajorgrids, xmajorgrids,
]
\fill[mycolor2!12,draw=none]
  (axis cs:-0.5,\yMinNorm) rectangle (axis cs:\exploreEndCEC+0.5,\yMaxNorm);
\addplot [color=mycolor2, line width=1.5pt] table[row sep=crcr]{%
0	21.2132034355964\\
1	9.48077016516599\\
2	2.12729774236618\\
3	3.32561050065416\\
4	1.97555928985174\\
5	0.595271746585285\\
6	2.41941971921537\\
7	1.20447779614633\\
8	1.63116294272808\\
9	0.381005400722809\\
10	1.52236629181327\\
11	1.4462748233376\\
12	1.47392752712804\\
13	1.75918274177432\\
14	2.80244402850699\\
15	2.44116798166673\\
16	0.409543993308893\\
17	0.918187387946841\\
18	3.78854144103001\\
19	2.13189533771111\\
20	2.71728758388557\\
21	2.09261710500357\\
22	2.75888016323964\\
23	1.98790243541202\\
24	0.954920763133843\\
25	2.51281447130344\\
26	3.45305321619226\\
27	1.01913024570199\\
28	0.97319769908553\\
29	0.97154453282467\\
30	2.15300212096701\\
31	2.25090359231126\\
32	1.27105106844081\\
33	0.37606582218967\\
34	0.732724893608207\\
35	3.68414314517416\\
36	2.32272659600168\\
37	1.14632000319405\\
38	1.44059266559553\\
39	1.8753412123591\\
40	1.55487427574816\\
41	1.3299843256896\\
42	0.705339131886646\\
43	0.263297976771175\\
44	0.614023389128065\\
45	0.519902812056371\\
46	0.45431186993246\\
47	0.76025552710524\\
48	1.27434476323232\\
49	0.0888330542868185\\
50	0.282615818656043\\
51	0.149532986582365\\
52	0.694165656677659\\
53	0.233511223778411\\
54	0.872580010770043\\
55	1.21931107609875\\
56	0.355032496437438\\
57	0.538524118672014\\
58	0.436013701281677\\
59	0.791265698882071\\
60	0.496195792553395\\
61	0.545795766017909\\
62	0.571807494744842\\
63	0.0784418456339577\\
64	0.301973214573684\\
65	0.970508966891256\\
66	0.563740934320541\\
67	0.578426982172433\\
68	0.236646970493254\\
69	0.309314329605196\\
70	0.816641974434498\\
71	0.134675361227865\\
72	0.587152608705899\\
73	0.402545160954177\\
74	0.204303967374051\\
75	0.306885107194401\\
76	0.55841362253547\\
77	0.825811636088872\\
78	0.066032189131571\\
79	0.344628449661856\\
80	1.15329176024382\\
81	0.291696027356094\\
82	0.692055027024523\\
83	0.921233000550089\\
84	1.3284301302518\\
85	1.01818066806333\\
86	0.674435197605075\\
87	0.175859493385209\\
88	0.908209181653638\\
89	1.06141038978647\\
90	0.344437863898316\\
91	0.698246104701626\\
92	0.211178512665974\\
93	0.392927071708644\\
94	0.957740952164392\\
95	1.07843398013526\\
96	0.644121219959778\\
97	0.0947457388835482\\
98	0.277836270957696\\
99	0.614500010892317\\
100	0.548528100930073\\
};
\node[anchor=north east,font=\scriptsize,fill=white,fill opacity=0.75,text opacity=1,inner sep=1pt] at (axis cs:\xMaxAll,\yMaxNorm) {Regret rate minimization};
\end{axis}

\begin{axis}[
  vaxisbase,
  at={(0,\tikzplotheight+\vplotspacing)},
  axis x line=none,
  axis y line*=right,
  ymin=\yMinB, ymax=\yMaxB,
  ylabel={$B,\hat B_t$},
  ylabel style={color=black},
]
\addplot [color=black, dashed, line width=1pt, forget plot] table[row sep=crcr]{%
0	0.707106781186547\\
1	0.707106781186547\\
99	0.707106781186547\\
};
\addplot [color=black, dashed, line width=1pt, forget plot] table[row sep=crcr]{%
0	1.5\\
1	1.5\\
99	1.5\\
};
\addplot [color=mycolor2, line width=1.5pt, dashdotted, forget plot] table[row sep=crcr]{%
0	0\\
1	1.19283003016223\\
2	1.04891396010133\\
3	0.816886470953187\\
4	0.794099690887516\\
5	0.832295197425764\\
6	0.867366180940223\\
7	0.868272540629537\\
8	0.899920770428165\\
9	0.907834292913705\\
10	0.893913917032278\\
11	0.872822465846713\\
12	0.692150647045449\\
13	0.697050804926342\\
14	0.677122978313274\\
15	0.693444157833101\\
16	0.665818245457226\\
17	0.609796590012926\\
18	0.730565718072527\\
19	0.731022828780404\\
20	0.683892447463401\\
21	0.68049860858508\\
22	0.756265671874857\\
23	0.717617249026681\\
24	0.71935305482173\\
25	0.662437595866256\\
26	0.689095573779038\\
27	0.685811649315012\\
28	0.687107929517469\\
29	0.687922060375922\\
30	0.681466042048147\\
31	0.634503225475415\\
32	0.631194565321781\\
33	0.631339371955034\\
34	0.633562957043534\\
35	0.670661024800582\\
36	0.653428557734722\\
37	0.648856839632826\\
38	0.639361025284733\\
39	0.634102395761189\\
40	0.627830756967396\\
41	0.626484722150443\\
42	0.626183455833804\\
43	0.626464581926164\\
44	0.625944897584356\\
45	0.625969088635238\\
46	0.626327051789972\\
47	0.625745546741613\\
48	0.630857817289155\\
49	0.631769100934664\\
50	0.631776184507669\\
51	0.631803455280557\\
52	0.631607864341739\\
53	0.632194202423843\\
54	0.632298734300996\\
55	0.62855007433372\\
56	0.62855977954647\\
57	0.628216244997692\\
58	0.627710602375122\\
59	0.627442612885109\\
60	0.62725162560137\\
61	0.627252224009268\\
62	0.628036957665359\\
63	0.629404185034343\\
64	0.62927502525963\\
65	0.62904477909496\\
66	0.6379075328166\\
67	0.637468112894767\\
68	0.637952332382631\\
69	0.637758242352651\\
70	0.639020827025969\\
71	0.641975280326661\\
72	0.641984616789406\\
73	0.641717953547909\\
74	0.641642753662938\\
75	0.641522323235303\\
76	0.640618457770007\\
77	0.640748569204715\\
78	0.64264340367285\\
79	0.642845341349149\\
80	0.642468573214913\\
81	0.650816358222553\\
82	0.649477611360323\\
83	0.648824895589306\\
84	0.640584435790037\\
85	0.641160624560685\\
86	0.639545233776749\\
87	0.639680122718765\\
88	0.638093542364747\\
89	0.634886581114321\\
90	0.635824748371929\\
91	0.640152552580043\\
92	0.640476903040608\\
93	0.640350137884036\\
94	0.638217098548848\\
95	0.632300556995586\\
96	0.636212403942534\\
97	0.638525580903844\\
98	0.638482259552376\\
99	0.638373964794906\\
};
\addplot [color=mycolor2, line width=1.5pt, dashdotted, forget plot] table[row sep=crcr]{%
0	0\\
1	1.40781915155643\\
2	1.30065956910623\\
3	1.57292831706717\\
4	1.56771143424035\\
5	1.55282119603935\\
6	1.62633373318882\\
7	1.62734635914121\\
8	1.61876811022372\\
9	1.61105912006568\\
10	1.61612273371448\\
11	1.57437987725078\\
12	1.53630759331126\\
13	1.54432875339401\\
14	1.58409799730687\\
15	1.58646173440062\\
16	1.57793656838858\\
17	1.52690499509741\\
18	1.64404721336209\\
19	1.64353946587957\\
20	1.67345421276653\\
21	1.67292115729147\\
22	1.67299046228663\\
23	1.67917920705224\\
24	1.67732372233773\\
25	1.68580270087597\\
26	1.66922859940593\\
27	1.63541560093626\\
28	1.63136872187112\\
29	1.63802100078485\\
30	1.58524872311349\\
31	1.59769274393745\\
32	1.59085515669116\\
33	1.59204020678532\\
34	1.59341037134024\\
35	1.5317754296042\\
36	1.54815829506498\\
37	1.55356254184272\\
38	1.55411507733206\\
39	1.57976073575909\\
40	1.58726625457542\\
41	1.59136822117208\\
42	1.58963616051984\\
43	1.58979953750486\\
44	1.58973869158287\\
45	1.58967651320735\\
46	1.58831659552987\\
47	1.58681499330251\\
48	1.58511767505117\\
49	1.5849135029638\\
50	1.58490088923203\\
51	1.584911911645\\
52	1.58483341790887\\
53	1.58488368840409\\
54	1.5852938251416\\
55	1.5858894166816\\
56	1.58526950934493\\
57	1.58605103322626\\
58	1.58571492567364\\
59	1.58573172841855\\
60	1.58572913455272\\
61	1.58574239968697\\
62	1.58766496687585\\
63	1.58695844291702\\
64	1.58687725405813\\
65	1.58784781946248\\
66	1.58580602091943\\
67	1.58560297453475\\
68	1.5858748736684\\
69	1.58616252767528\\
70	1.58722271040659\\
71	1.58662184120171\\
72	1.5866199457014\\
73	1.58689623185887\\
74	1.5869381979195\\
75	1.58680555080941\\
76	1.58669392828036\\
77	1.58678997930322\\
78	1.58601596640565\\
79	1.58613607849988\\
80	1.588574956714\\
81	1.58243255283976\\
82	1.58093842938693\\
83	1.58457914868457\\
84	1.58241701101864\\
85	1.58298538469635\\
86	1.57982701758261\\
87	1.58015082708234\\
88	1.57950615154957\\
89	1.5805119646022\\
90	1.57881332472596\\
91	1.5794089619837\\
92	1.57988156593922\\
93	1.58086239954918\\
94	1.57671814152224\\
95	1.57997682502464\\
96	1.5759543546562\\
97	1.57477305942488\\
98	1.57494062397541\\
99	1.57729340724471\\
};
\end{axis}

\begin{axis}[
  vaxisbase,
  at={(0,0)},
  axis x line*=bottom,
  axis y line*=left,
  xlabel style={font=\color{mycolor4}},
  xlabel={$t$},
  ymin=\yMinNorm, ymax=\yMaxNorm,
  ylabel={$\|x(t)\|$},
  ymajorgrids, xmajorgrids,
]
\fill[mycolor1!15,draw=none]
  (axis cs:-0.5,\yMinNorm) rectangle (axis cs:\exploreEndMinimax+0.5,\yMaxNorm);
\addplot [color=mycolor1, line width=1.5pt] table[row sep=crcr]{%
0	21.2132034355964\\
1	16.3904764987823\\
2	5.16529706561991\\
3	0.924548781072723\\
4	0.284404497684786\\
5	0.794649760851908\\
6	0.587763346537719\\
7	0.607469867607245\\
8	0.722258431804461\\
9	0.776631049097467\\
10	0.338769667432587\\
11	0.34848619921505\\
12	1.2803163205152\\
13	0.159962527834816\\
14	0.303576264207952\\
15	0.541140433537663\\
16	0.698431775519155\\
17	1.28348234330741\\
18	1.16113717727528\\
19	0.731094400706057\\
20	0.472856183664833\\
21	0.13236844577125\\
22	1.08627794556096\\
23	0.697524460224882\\
24	0.103654228230196\\
25	0.78452597406824\\
26	0.55587977513952\\
27	0.610760483647152\\
28	1.4087943230875\\
29	0.982242481180592\\
30	0.756398773601266\\
31	0.561871684726533\\
32	0.439107572957354\\
33	0.273685039867589\\
34	0.120866164452268\\
35	0.54652528118669\\
36	1.24365740419087\\
37	0.452133421083297\\
38	1.16700230186282\\
39	1.13182643894164\\
40	0.315375134539945\\
41	0.943484882308168\\
42	0.587636501476821\\
43	0.294597297983654\\
44	0.60606393782726\\
45	0.526496784551613\\
46	0.405178723682555\\
47	0.710573732032021\\
48	1.29850101824664\\
49	0.135272454125047\\
50	0.296771103940989\\
51	0.153289453319766\\
52	0.702264599184357\\
53	0.256543864338935\\
54	0.882874296770518\\
55	1.20925335774531\\
56	0.335927937634194\\
57	0.527365460397709\\
58	0.426609530678709\\
59	0.76398739650487\\
60	0.490146871497777\\
61	0.54693359455989\\
62	0.615211220011897\\
63	0.0411270864324796\\
64	0.270930218052249\\
65	0.984422815191753\\
66	0.585605235019988\\
67	0.530236355564285\\
68	0.256099588761269\\
69	0.329076756335158\\
70	0.832766993110485\\
71	0.160334692270455\\
72	0.628031838684347\\
73	0.416368201880387\\
74	0.212569535293545\\
75	0.309454763551173\\
76	0.540801822109698\\
77	0.832572915700224\\
78	0.058114746225178\\
79	0.36638542935879\\
80	1.16724096920471\\
81	0.173248400074676\\
82	0.62839409901681\\
83	0.950875562462745\\
84	1.31085719173651\\
85	1.04615474144187\\
86	0.612206268716933\\
87	0.161680858781913\\
88	0.897659222519342\\
89	1.04704616307003\\
90	0.307818016446733\\
91	0.740986424089304\\
92	0.235853738854503\\
93	0.419829440022134\\
94	0.925877853347038\\
95	1.05098813124897\\
96	0.587008186712185\\
97	0.0264291680499859\\
98	0.263943770530279\\
99	0.633640760982663\\
100	0.508098825718515\\
};
\node[anchor=north east,font=\scriptsize,fill=white,fill opacity=0.75,text opacity=1,inner sep=1pt] at (axis cs:\xMaxAll,\yMaxNorm) {Minimax approach};
\end{axis}

\begin{axis}[
  vaxisbase,
  at={(0,0)},
  axis x line=none,
  axis y line*=right,
  ymin=\yMinB, ymax=\yMaxB,
  ylabel={$B,\hat B_t$},
  ylabel style={color=black},
]
\addplot [color=black, dashed, line width=1pt, forget plot] table[row sep=crcr]{%
0	0.707106781186547\\
1	0.707106781186547\\
99	0.707106781186547\\
};
\addplot [color=black, dashed, line width=1pt, forget plot] table[row sep=crcr]{%
0	1.5\\
1	1.5\\
99	1.5\\
};
\addplot [color=mycolor1, line width=1.5pt, dashdotted, forget plot] table[row sep=crcr]{%
0	2.99999290563833\\
1	0.614256033408615\\
2	0.627985194898341\\
3	0.629207116413599\\
4	0.629210213171436\\
5	0.629011119686205\\
6	0.628539505386905\\
7	0.629395761507371\\
8	0.6281797688441\\
9	0.627790323787212\\
10	0.62769452440029\\
11	0.627690129673219\\
12	0.625403325515599\\
13	0.625443300654972\\
14	0.625603460971687\\
15	0.624885275628483\\
16	0.624951726722555\\
17	0.626192607643479\\
18	0.631043704193958\\
19	0.63140218631423\\
20	0.632486123324143\\
21	0.632864927356321\\
22	0.633239053650384\\
23	0.632517505059923\\
24	0.633004085526091\\
25	0.632467761922154\\
26	0.631471061334314\\
27	0.631529118726043\\
28	0.629692120311447\\
29	0.630318200228506\\
30	0.631087368251531\\
31	0.637109674100415\\
32	0.637428495134051\\
33	0.637407631601627\\
34	0.637462538929123\\
35	0.637189208699343\\
36	0.633937092638592\\
37	0.638328997084456\\
38	0.634190196403034\\
39	0.633665879667927\\
40	0.638574540223895\\
41	0.637899371896958\\
42	0.637372334940108\\
43	0.637791206515065\\
44	0.636632080373289\\
45	0.636624636279014\\
46	0.637439356678477\\
47	0.636279315076869\\
48	0.647123516362507\\
49	0.648959303210631\\
50	0.648915003902308\\
51	0.648984428036431\\
52	0.64831931416552\\
53	0.649575619613741\\
54	0.649807233731602\\
55	0.64139074915755\\
56	0.641430870348023\\
57	0.640764543044987\\
58	0.639614109742196\\
59	0.639536704600143\\
60	0.639159493612056\\
61	0.639131851341377\\
62	0.640887569896982\\
63	0.644257849377802\\
64	0.644227505908363\\
65	0.64358453478199\\
66	0.663650767130104\\
67	0.663485766067828\\
68	0.664640321213038\\
69	0.664152278240805\\
70	0.667273829945333\\
71	0.674205700402626\\
72	0.673411096512418\\
73	0.672718945453421\\
74	0.672548041261635\\
75	0.672144152924382\\
76	0.670262677920027\\
77	0.670261480061446\\
78	0.674505473906457\\
79	0.674572700338208\\
80	0.673568938529681\\
81	0.692573540738784\\
82	0.690725069634084\\
83	0.689133452402644\\
84	0.669861534788997\\
85	0.669087323604684\\
86	0.665660657800341\\
87	0.665858312006644\\
88	0.662257472699342\\
89	0.655520371756998\\
90	0.657768800748438\\
91	0.666374788087773\\
92	0.667193188987252\\
93	0.666861221571221\\
94	0.661685436780807\\
95	0.649333316243923\\
96	0.659246623559203\\
97	0.663914797990063\\
98	0.663901764658952\\
99	0.663598502055939\\
};
\addplot [color=mycolor1, line width=1.5pt, dashdotted, forget plot] table[row sep=crcr]{%
0	0\\
1	1.51810724012342\\
2	1.50327983597957\\
3	1.49923846154449\\
4	1.49923905220079\\
5	1.49931016789396\\
6	1.49868418218846\\
7	1.49964021864125\\
8	1.49957713649618\\
9	1.49993196944386\\
10	1.4995391267087\\
11	1.49954053420878\\
12	1.4991467902697\\
13	1.4992331010772\\
14	1.49887937011169\\
15	1.4987000047986\\
16	1.49871801607015\\
17	1.49981211920717\\
18	1.50572177676472\\
19	1.50609200797781\\
20	1.50487932439277\\
21	1.50469655514305\\
22	1.50475320935088\\
23	1.50511762244989\\
24	1.50467118684217\\
25	1.50496814415622\\
26	1.504997133544\\
27	1.50552363372032\\
28	1.51119983712893\\
29	1.51546932858209\\
30	1.51871308800704\\
31	1.51663489699132\\
32	1.51719493572242\\
33	1.51642977330445\\
34	1.5165114195092\\
35	1.51702838817381\\
36	1.52011714538702\\
37	1.51503498248653\\
38	1.51532391725972\\
39	1.51728472922069\\
40	1.51164562885883\\
41	1.51355665032691\\
42	1.51020899519627\\
43	1.51045658372847\\
44	1.51032288101906\\
45	1.51030380454597\\
46	1.50729293223807\\
47	1.50433963130549\\
48	1.50077520023754\\
49	1.50038145292364\\
50	1.50045325294309\\
51	1.5004845075728\\
52	1.50022127457224\\
53	1.50033424545138\\
54	1.50153247835477\\
55	1.50295784695342\\
56	1.50112235889872\\
57	1.50263629711561\\
58	1.50187106788783\\
59	1.50187604969716\\
60	1.50187122882064\\
61	1.50194692195736\\
62	1.50627586023652\\
63	1.50459593685386\\
64	1.50457922840836\\
65	1.50693700001875\\
66	1.50243675640626\\
67	1.50236570364692\\
68	1.50302125130942\\
69	1.5037514454239\\
70	1.50643481723259\\
71	1.50507966007158\\
72	1.50524142063911\\
73	1.50595672697752\\
74	1.50606597474651\\
75	1.50562759763691\\
76	1.50540022497425\\
77	1.50544055384821\\
78	1.5037545871976\\
79	1.50380041884169\\
80	1.50916544303939\\
81	1.49546939038207\\
82	1.49345667812268\\
83	1.5017849078845\\
84	1.49683924613304\\
85	1.49619781179858\\
86	1.4894674153618\\
87	1.48995836477026\\
88	1.48850998560266\\
89	1.49064795885438\\
90	1.48662456275341\\
91	1.48783908607487\\
92	1.48904863771624\\
93	1.49151515088989\\
94	1.48165741687275\\
95	1.488658248963\\
96	1.47856154749803\\
97	1.47629271545453\\
98	1.47629468383664\\
99	1.48147589656048\\
};
\end{axis}

\end{tikzpicture}%

%% file: figures/Ljung-colored_auto.tex
\definecolor{mycolor1}{rgb}{0.00000,0.44700,0.74100}
\definecolor{mycolor2}{rgb}{0.85000,0.32500,0.09800}
\definecolor{mycolor3}{rgb}{0.92900,0.69400,0.12500}
\definecolor{mycolor4}{rgb}{0.12941,0.12941,0.12941}%

\setlength{\tikzplotheight}{1.5cm}
\setlength{\vplotspacing}{.25cm}
\setlength{\tikzplotwidth}{.75\linewidth}

\def\exploreEndCEC{33}
\def\exploreEndWRLS{33}
\def\exploreEndMinimax{0}

\def\xMinAll{0}
\def\xMaxAll{60}
\def\yMinNorm{0}
\def\yMaxNorm{30}
\def\yMinB{-1}
\def\yMaxB{4}

\begin{tikzpicture}

\pgfplotsset{
  vaxisbase/.style={
    width=\tikzplotwidth,
    height=\tikzplotheight,
    scale only axis,
    xmin=\xMinAll, xmax=\xMaxAll,
    xtick={0,20,40,60,80,100,120,140},
    axis background/.style={fill=none},
    clip=true,
  },
}

\begin{axis}[
  vaxisbase,
  at={(0,2\tikzplotheight+2\vplotspacing)},
  axis x line*=bottom,
  axis y line*=left,
  xticklabels={},
  ymin=\yMinNorm, ymax=\yMaxNorm,
  ylabel={$\|x(t)\|$},
  legend style={legend cell align=left, align=left, font=\scriptsize,
    at={(0.05,0.87)}, anchor=north west, draw=none, fill=none, /tikz/every even column/.append style={column sep=4pt}},
  ymajorgrids, xmajorgrids,
]
\fill[mycolor3!15,draw=none]
  (axis cs:-0.5,\yMinNorm) rectangle (axis cs:\exploreEndWRLS+0.5,\yMaxNorm);
\addplot [color=mycolor3, line width=1.5pt] table[row sep=crcr]{%
0	21.2132034355964\\
1	64.6095394755015\\
2	57.2637131801838\\
3	6.75906696317982\\
4	2.48055822586123\\
5	2.18329298801473\\
6	2.75947813992564\\
7	0.723391751809531\\
8	0.941319641471983\\
9	0.556443470312115\\
10	0.958048643110962\\
11	1.20545529470838\\
12	1.50069635185432\\
13	3.28721146235242\\
14	3.48394750240789\\
15	3.80727567779889\\
16	3.78879444463284\\
17	6.16150584771515\\
18	4.11826969525066\\
19	1.00781212534867\\
20	0.984283665770928\\
21	4.51670832955968\\
22	5.56378934719391\\
23	3.77380580873368\\
24	0.972444879472769\\
25	3.27745805062315\\
26	3.34968256762534\\
27	3.94535303380731\\
28	3.7671107811836\\
29	0.390198286674626\\
30	2.6663032233124\\
31	4.72007904161375\\
32	2.66168159055624\\
33	2.09059336519472\\
34	3.19978817640529\\
35	1.1623083977634\\
36	3.50937235178631\\
37	3.34442631318638\\
38	1.82488515227169\\
39	2.3530130006339\\
40	3.75255706308658\\
41	1.49496027485681\\
42	2.19829871050416\\
43	3.8488080528117\\
44	2.27745102956011\\
45	3.14200173873563\\
46	1.89186912842946\\
47	1.98786981828637\\
48	2.1364185121\\
49	1.7722878597331\\
50	2.41155013465209\\
51	0.799823996936624\\
52	1.10077939052523\\
53	0.558954739722421\\
54	0.319414629767708\\
55	1.6616395337434\\
56	2.28183341038609\\
57	1.05878744458169\\
58	0.472530242076819\\
59	3.99012186187493\\
60	4.76406428689985\\
61	0.89495683713635\\
62	1.50536402315726\\
63	1.85949561013294\\
64	1.71394206821719\\
65	1.29963703753296\\
66	1.12701694244471\\
67	2.10584445757074\\
68	3.01173944851137\\
69	0.977584713026481\\
70	0.487635210980954\\
71	0.318936643189354\\
72	0.965192522139988\\
73	0.998854531822168\\
74	4.03951536109377\\
75	5.01001474506385\\
76	3.99933229419661\\
77	1.29057481243539\\
78	1.73849501115738\\
79	2.34360222007174\\
80	1.99453912944011\\
81	1.15677182059251\\
82	0.719595959520605\\
83	3.08031338330858\\
84	3.67927660587126\\
85	2.23729572607572\\
86	3.75445397599701\\
87	7.14772360069532\\
88	4.8523068331243\\
89	2.33427532183046\\
90	3.30075995820799\\
91	1.46390836639688\\
92	2.39252304845471\\
93	2.73483233128031\\
94	3.21311774260376\\
95	3.27179945254474\\
96	3.26804506793843\\
97	1.07863484447695\\
98	1.57426748387123\\
99	2.80643048637549\\
100	2.98476321359719\\
};
\addlegendentry{$\|x(t)\|$}
\node[anchor=north east,font=\scriptsize,fill=white,fill opacity=0.75,text opacity=1,inner sep=1pt] at (axis cs:\xMaxAll,\yMaxNorm) {Weighted recursive least squares};
\end{axis}

\begin{axis}[
  vaxisbase,
  at={(0,2\tikzplotheight+2\vplotspacing)},
  axis x line=none,
  axis y line*=right,
  ymin=\yMinB, ymax=\yMaxB,
  ylabel={$B,\hat B_t$},
  ylabel style={color=black},
  legend style={legend cell align=left, align=left, font=\scriptsize,
    at={(1.02,0.87)}, anchor=north east, draw=none, fill=none, /tikz/every even column/.append style={column sep=4pt}},
  legend columns=2,
]
\addplot [color=black, dashed, line width=1pt, forget plot] table[row sep=crcr]{%
0	1.95870836561711\\
1	1.95870836561711\\
99	1.95870836561711\\
};
\addplot [color=black, dashed, line width=1pt] table[row sep=crcr]{%
0	0\\
1	0\\
99	0\\
};
\addlegendentry{$B$}
\addplot [color=mycolor3, line width=1.5pt, dashdotted] table[row sep=crcr]{%
0	0.00761440008939498\\
1	1.87270948393949\\
2	1.86987510683528\\
3	1.84672120742473\\
4	1.84770901218146\\
5	1.84800793586921\\
6	1.84807075693247\\
7	1.84837407000974\\
8	1.84837723329386\\
9	1.84843825870398\\
10	1.84852395851053\\
11	1.8484595010924\\
12	1.8484325228759\\
13	1.84879896896366\\
14	1.84880235752015\\
15	1.84879286196312\\
16	1.84878526992369\\
17	1.8503149858244\\
18	1.84995135167518\\
19	1.84936896720147\\
20	1.84936866683949\\
21	1.85006029246204\\
22	1.85005985816289\\
23	1.85050802295306\\
24	1.84931882532506\\
25	1.84932850835304\\
26	1.84938133795256\\
27	1.84935511183936\\
28	1.8493744046129\\
29	1.84975108995732\\
30	1.84982852615105\\
31	1.85033875849997\\
32	1.85102945158899\\
33	1.85133137114874\\
34	1.85131743480536\\
35	1.8512421813599\\
36	1.85161331527597\\
37	1.85178050167435\\
38	1.85179966203502\\
39	1.85179578384728\\
40	1.85170432469627\\
41	1.85165234984655\\
42	1.8517048611696\\
43	1.85187069059627\\
44	1.85178779331472\\
45	1.85190077055015\\
46	1.85227989480446\\
47	1.85276376868105\\
48	1.85266503463373\\
49	1.85285668151728\\
50	1.85287058643735\\
51	1.85290548199225\\
52	1.8529140469015\\
53	1.85294294017965\\
54	1.85293478164725\\
55	1.85281439801246\\
56	1.85293172924815\\
57	1.85313983265271\\
58	1.85317418522819\\
59	1.85353716043145\\
60	1.85346965232703\\
61	1.85381230855079\\
62	1.85383341797118\\
63	1.85384313191226\\
64	1.85384684817584\\
65	1.85422875077822\\
66	1.85413768973525\\
67	1.85420914866581\\
68	1.85415247094386\\
69	1.85384503707632\\
70	1.85381961751504\\
71	1.85378834511312\\
72	1.85380953104774\\
73	1.85372398150514\\
74	1.85419845325282\\
75	1.85427105960067\\
76	1.85439650637159\\
77	1.85450093902537\\
78	1.8545168822136\\
79	1.85447292616197\\
80	1.85474118686344\\
81	1.85464538535321\\
82	1.85452027414954\\
83	1.85504662872338\\
84	1.85501556754224\\
85	1.85506739037327\\
86	1.85527759369297\\
87	1.85676972671938\\
88	1.8574657069833\\
89	1.85744459540869\\
90	1.85764721349508\\
91	1.85821351672598\\
92	1.85829076460247\\
93	1.85829586163901\\
94	1.85797969161764\\
95	1.85783634602634\\
96	1.8577740476444\\
97	1.85817248684169\\
98	1.85792825831414\\
99	1.85819261778959\\
};
\addlegendentry{$\hat{B}_t$}
\addplot [color=mycolor3, line width=1.5pt, dashdotted, forget plot] table[row sep=crcr]{%
0	0.695338824921826\\
1	0.150586395564603\\
2	0.113430639579103\\
3	0.101065906822133\\
4	0.100651365813905\\
5	0.100628243794538\\
6	0.10063204744564\\
7	0.100633209732174\\
8	0.100647987262951\\
9	0.10062991974456\\
10	0.100683860224617\\
11	0.100657070965437\\
12	0.100662521886139\\
13	0.100702346718159\\
14	0.100703299483481\\
15	0.100705750912409\\
16	0.100725679256468\\
17	0.100569302135579\\
18	0.100510424515935\\
19	0.100041753944943\\
20	0.100031714370465\\
21	0.100165912728746\\
22	0.100134457344103\\
23	0.100057047438436\\
24	0.0995364009443153\\
25	0.0993492903756397\\
26	0.0993442690585113\\
27	0.0993025180767555\\
28	0.099307963283463\\
29	0.099273615593568\\
30	0.0991954031313349\\
31	0.099130673459134\\
32	0.0990580408675801\\
33	0.0991123971570315\\
34	0.0989908473707076\\
35	0.0989930776072019\\
36	0.0990748171708019\\
37	0.0989882253466204\\
38	0.0989979050767668\\
39	0.0990120456256548\\
40	0.0989003446705765\\
41	0.0988319268001213\\
42	0.0988539454993229\\
43	0.0988081076566516\\
44	0.0987556993383366\\
45	0.0987744518105487\\
46	0.0987450739063823\\
47	0.0988046560235192\\
48	0.0988078773755978\\
49	0.0988929391481406\\
50	0.0988812182877353\\
51	0.0988962360475963\\
52	0.0989008130432144\\
53	0.098878522056408\\
54	0.0988774761981902\\
55	0.0988305840208284\\
56	0.0987448991959855\\
57	0.0987747104458046\\
58	0.0987635743212686\\
59	0.0987929939189962\\
60	0.0986768393310029\\
61	0.0986144748991414\\
62	0.0986347331041711\\
63	0.0986211939564349\\
64	0.0986214729134271\\
65	0.0986927559873085\\
66	0.0986896660372003\\
67	0.0987319576106238\\
68	0.0987054886055188\\
69	0.0986018644078359\\
70	0.0985840476008149\\
71	0.098546905808003\\
72	0.0985491435557447\\
73	0.098549437781814\\
74	0.0986952893475667\\
75	0.0986808021173982\\
76	0.0986790208964221\\
77	0.0986064760354352\\
78	0.0986031585439185\\
79	0.0985826035187012\\
80	0.0984702955827197\\
81	0.098501834678979\\
82	0.0983663816946003\\
83	0.0984156308385995\\
84	0.0983690690676215\\
85	0.0980517228526329\\
86	0.0980253972624108\\
87	0.0981861836622333\\
88	0.0980506153825242\\
89	0.0979839402635282\\
90	0.0979534016624323\\
91	0.0980005677070289\\
92	0.0979845028150087\\
93	0.0979800043807686\\
94	0.0979086729306896\\
95	0.0978808876006801\\
96	0.0978785473589544\\
97	0.0979073491897702\\
98	0.0978827152286698\\
99	0.0979574575310481\\
};
\end{axis}

\begin{axis}[
  vaxisbase,
  at={(0,\tikzplotheight+\vplotspacing)},
  axis x line*=bottom,
  axis y line*=left,
  xticklabels={},
  ymin=\yMinNorm, ymax=\yMaxNorm,
  ylabel={$\|x(t)\|$},
  ymajorgrids, xmajorgrids,
]
\fill[mycolor2!12,draw=none]
  (axis cs:-0.5,\yMinNorm) rectangle (axis cs:\exploreEndCEC+0.5,\yMaxNorm);
\addplot [color=mycolor2, line width=1.5pt] table[row sep=crcr]{%
0	21.2132034355964\\
1	23.840988323379\\
2	19.5426046709935\\
3	23.3766444504199\\
4	26.4817048853556\\
5	23.5077584265034\\
6	25.2261455532626\\
7	29.2011034234344\\
8	20.3387603822775\\
9	21.4256561897699\\
10	22.5816246063992\\
11	17.5772793330782\\
12	20.9742349745695\\
13	21.2752782471633\\
14	16.8464375759029\\
15	19.0972812823426\\
16	17.7559267341251\\
17	17.0584186296126\\
18	19.2436348752675\\
19	17.4237632083838\\
20	15.648447022888\\
21	16.7414519574525\\
22	11.0745653038883\\
23	11.8089250986234\\
24	12.550385592086\\
25	9.42609125500035\\
26	11.7789612065416\\
27	12.1709746636097\\
28	7.97937614907511\\
29	7.72109765921144\\
30	9.04799260723141\\
31	7.61608104423812\\
32	6.87847845325347\\
33	7.3842916685956\\
34	6.12696736142609\\
35	3.37371536451648\\
36	0.842728209771667\\
37	0.232093421836301\\
38	0.379736048826414\\
39	0.306591173355064\\
40	0.297406620765647\\
41	0.476718826946102\\
42	0.355623783995084\\
43	0.376444920008164\\
44	0.23855265486306\\
45	0.267703242054937\\
46	0.351249016920388\\
47	0.440797337814046\\
48	0.671998304643044\\
49	0.490488095161776\\
50	0.475123635633824\\
51	0.137450207221379\\
52	0.207880822014721\\
53	0.114788502501017\\
54	0.227991925562662\\
55	0.369116445471081\\
56	0.161878499903383\\
57	0.178129785505394\\
58	0.263873056442573\\
59	0.104965023697618\\
60	0.29729553567228\\
61	0.291359977322434\\
62	0.107084657176289\\
63	0.215435760812113\\
64	0.252629646915879\\
65	0.485116066320915\\
66	0.693655549965737\\
67	0.322004390627829\\
68	0.349439256109414\\
69	0.535118246437216\\
70	0.32252065844395\\
71	0.384024224142484\\
72	0.504005814202302\\
73	0.577599576364835\\
74	0.440556691755106\\
75	0.428779929083814\\
76	0.452336265228583\\
77	0.186946753737045\\
78	0.274926223407326\\
79	0.408171930138436\\
80	0.23554748097049\\
81	0.231365699144197\\
82	0.35759790778972\\
83	0.595162558258695\\
84	0.405128140339949\\
85	0.304307047338944\\
86	0.477100912088051\\
87	0.580771686677661\\
88	0.577274573320418\\
89	0.176859905197388\\
90	0.311323670601865\\
91	0.321747771171949\\
92	0.320061964392554\\
93	0.421705857164843\\
94	0.678570168789669\\
95	0.580876195408557\\
96	0.381416044990931\\
97	0.303432324109357\\
98	0.565189075683303\\
99	0.7983562151542\\
100	0.859323502583408\\
};
\node[anchor=north east,font=\scriptsize,fill=white,fill opacity=0.75,text opacity=1,inner sep=1pt] at (axis cs:\xMaxAll,\yMaxNorm) {Regret rate minimization};
\end{axis}

\begin{axis}[
  vaxisbase,
  at={(0,\tikzplotheight+\vplotspacing)},
  axis x line=none,
  axis y line*=right,
  ymin=\yMinB, ymax=\yMaxB,
  ylabel={$B,\hat B_t$},
  ylabel style={color=black},
]
\addplot [color=black, dashed, line width=1pt, forget plot] table[row sep=crcr]{%
0	1.95870836561711\\
1	1.95870836561711\\
99	1.95870836561711\\
};
\addplot [color=black, dashed, line width=1pt, forget plot] table[row sep=crcr]{%
0	0\\
1	0\\
99	0\\
};
\addplot [color=mycolor2, line width=1.5pt, dashdotted, forget plot] table[row sep=crcr]{%
0	0\\
1	1.78203285645697\\
2	1.79315773828539\\
3	2.08636220002668\\
4	2.04745240027094\\
5	2.10171697576188\\
6	2.05877526959949\\
7	1.9989774290245\\
8	1.97924474066518\\
9	1.99056402649783\\
10	1.98828836549013\\
11	1.98387800667696\\
12	1.99553803560749\\
13	1.99819335895758\\
14	1.99842675611176\\
15	1.91499422527727\\
16	1.93648442086891\\
17	1.93843882392007\\
18	1.91814264477228\\
19	1.89068367975443\\
20	1.91635601785053\\
21	1.91766326373239\\
22	1.93571416621944\\
23	1.93778269102757\\
24	1.93403930542195\\
25	1.93141373910267\\
26	1.93611556922084\\
27	1.94703068152102\\
28	1.94358450547916\\
29	1.94419012443473\\
30	1.94456435866805\\
31	1.94478050053576\\
32	1.94935449401554\\
33	1.95178083931016\\
34	1.95200362532766\\
35	1.93665779925776\\
36	1.92997456089596\\
37	1.93054560696478\\
38	1.931588060167\\
39	1.93160176744293\\
40	1.93063981880669\\
41	1.93091002295282\\
42	1.93326144600594\\
43	1.93314921285938\\
44	1.93216563197344\\
45	1.93279714178543\\
46	1.93396575150845\\
47	1.93624480837541\\
48	1.93941805106763\\
49	1.94352939258481\\
50	1.94440514305256\\
51	1.94510538335509\\
52	1.94510477124529\\
53	1.94512641633825\\
54	1.94542600657605\\
55	1.94573132904056\\
56	1.94600379661903\\
57	1.9457031366754\\
58	1.94561471511181\\
59	1.94585853961197\\
60	1.94591195968804\\
61	1.94595831059996\\
62	1.94584783196624\\
63	1.94588079257056\\
64	1.94563204260551\\
65	1.94782362484262\\
66	1.94907251804445\\
67	1.95335111976956\\
68	1.9518792937291\\
69	1.95055989004476\\
70	1.95164748621235\\
71	1.94922693630012\\
72	1.94937090312689\\
73	1.95702224897066\\
74	1.95889469923294\\
75	1.95580162026461\\
76	1.95632968556513\\
77	1.95658513212947\\
78	1.95567600729448\\
79	1.95536422213034\\
80	1.95557378717226\\
81	1.9570003138199\\
82	1.95899288439677\\
83	1.9576943593648\\
84	1.95976954206079\\
85	1.95854803749895\\
86	1.96166425668496\\
87	1.9708928341165\\
88	1.97049660847761\\
89	1.96802131191867\\
90	1.96878456978328\\
91	1.97233906363937\\
92	1.97198854164216\\
93	1.96810875529815\\
94	1.97523404411677\\
95	1.98344198224818\\
96	1.98870565778899\\
97	1.9930042325965\\
98	1.99797503832539\\
99	2.01480562211384\\
};
\addplot [color=mycolor2, line width=1.5pt, dashdotted, forget plot] table[row sep=crcr]{%
0	0\\
1	0\\
2	-0.0291401858848765\\
3	-0.0587613828016277\\
4	-0.112689540558811\\
5	-0.0531534546399648\\
6	-0.0503325186298949\\
7	-0.0547837310130086\\
8	-0.0616429172997634\\
9	-0.0633511814314356\\
10	-0.0662657225649375\\
11	-0.0674749366148625\\
12	-0.0615460051466686\\
13	-0.0590580259179978\\
14	-0.0591307195700895\\
15	-0.0417185940730136\\
16	-0.0284661215656823\\
17	-0.0276852884452642\\
18	-0.0206784647498313\\
19	-0.0445806053408727\\
20	-0.0381391474204999\\
21	-0.0341166185613215\\
22	-0.0247383797633555\\
23	-0.0254710012763745\\
24	-0.0274209438491399\\
25	-0.0307690262277028\\
26	-0.0252207524953875\\
27	-0.020314934336748\\
28	-0.0229059455949326\\
29	-0.0228153472076855\\
30	-0.0224200944048287\\
31	-0.0223052964027659\\
32	-0.0217990635099207\\
33	-0.0132509587644972\\
34	-0.0129923684156557\\
35	-0.00468378393288641\\
36	-0.0129412701467277\\
37	-0.0137977531903027\\
38	-0.0134356551463268\\
39	-0.0129789242476912\\
40	-0.0137430872783728\\
41	-0.0137474159939427\\
42	-0.0121545996201856\\
43	-0.0120159439676653\\
44	-0.0123032705156836\\
45	-0.0119051778074813\\
46	-0.0117719545182561\\
47	-0.010606851466478\\
48	-0.00985980836585932\\
49	-0.00641420161580583\\
50	-0.00662955936792143\\
51	-0.00570428875505064\\
52	-0.00570244175444299\\
53	-0.00544737185228462\\
54	-0.00532404820553042\\
55	-0.00522689049424397\\
56	-0.00424918161230136\\
57	-0.00462824224020092\\
58	-0.00458917829880075\\
59	-0.00373962332181\\
60	-0.00370452112587401\\
61	-0.00371300133039911\\
62	-0.00416799520004093\\
63	-0.00452231036115888\\
64	-0.00387501318720096\\
65	-0.00256403261504979\\
66	-0.00227608867249065\\
67	0.00147097561198039\\
68	0.00124313458426532\\
69	0.000840292167303339\\
70	0.00392308554686384\\
71	0.00228324397814824\\
72	0.00230453275491912\\
73	0.00442217156412188\\
74	0.00836662261500557\\
75	0.00561649269073551\\
76	0.00568512287024742\\
77	0.00479207847714388\\
78	0.00377562302590572\\
79	0.00365382192412738\\
80	0.00322552158888297\\
81	0.00246730372087071\\
82	0.00499697158760109\\
83	0.00408525283809733\\
84	0.00390186784357854\\
85	0.00270702309319015\\
86	0.00220465950418457\\
87	0.00769703148476898\\
88	0.00927421544358017\\
89	0.00682734866145283\\
90	0.00665871674428966\\
91	0.0087535578971962\\
92	0.009252266463334\\
93	0.00823160326337832\\
94	0.0113000299642603\\
95	0.0167073701078553\\
96	0.0175750950448231\\
97	0.0199266638221026\\
98	0.0206689006572617\\
99	0.0280530711152676\\
};
\end{axis}

\begin{axis}[
  vaxisbase,
  at={(0,0)},
  axis x line*=bottom,
  axis y line*=left,
  xlabel style={font=\color{mycolor4}},
  xlabel={$t$},
  ymin=\yMinNorm, ymax=\yMaxNorm,
  ylabel={$\|x(t)\|$},
  ymajorgrids, xmajorgrids,
]
\fill[mycolor1!15,draw=none]
  (axis cs:-0.5,\yMinNorm) rectangle (axis cs:\exploreEndMinimax+0.5,\yMaxNorm);
\addplot [color=mycolor1, line width=1.5pt] table[row sep=crcr]{%
0	21.2132034355964\\
1	25.1495702468866\\
2	23.8651289377203\\
3	3.06727375972734\\
4	0.661466218872192\\
5	0.779541512732488\\
6	0.801936997541404\\
7	0.485030172362022\\
8	0.422991990119919\\
9	0.466829844386456\\
10	0.353307361446705\\
11	0.57213462098653\\
12	0.834123654335157\\
13	0.856783553101923\\
14	0.538020396166058\\
15	0.483182195471361\\
16	0.779275374565437\\
17	0.883673632029094\\
18	0.358107523041645\\
19	0.778825850356191\\
20	0.969902516518548\\
21	0.2021309494781\\
22	0.519260047533005\\
23	0.403804746882024\\
24	1.06288527033569\\
25	1.34744267180789\\
26	0.668998045402533\\
27	0.158587137631725\\
28	0.209583515270065\\
29	0.172730825472362\\
30	0.142734470013913\\
31	0.169366341087118\\
32	0.250225513683825\\
33	0.212268088899932\\
34	0.336634840781523\\
35	0.400790530251575\\
36	0.334419736901281\\
37	0.523210351490395\\
38	0.397984820713347\\
39	0.395202739647256\\
40	0.37219083182824\\
41	0.537595648334225\\
42	0.493972285918676\\
43	0.518582544658702\\
44	0.212188658477787\\
45	0.349908521893997\\
46	0.495153865661201\\
47	0.588239774466587\\
48	0.904386830853706\\
49	0.798935769681921\\
50	0.768917182619814\\
51	0.199802014591958\\
52	0.257905526068059\\
53	0.0610419600231327\\
54	0.214334726570961\\
55	0.430067796795746\\
56	0.322646919456979\\
57	0.264374468488297\\
58	0.263883277725949\\
59	0.114318114513779\\
60	0.297706685165753\\
61	0.241734091535368\\
62	0.11044577347422\\
63	0.286861700497721\\
64	0.221298043088202\\
65	0.478309667707875\\
66	0.82905640275908\\
67	0.597546665769957\\
68	0.467399965388644\\
69	0.534987244177554\\
70	0.52647830967423\\
71	0.483616507756594\\
72	0.534685776301589\\
73	0.726932950372366\\
74	0.676740761710843\\
75	0.609583965346562\\
76	0.347476607351511\\
77	0.0928668869578304\\
78	0.323542868750771\\
79	0.391283371725402\\
80	0.14341510096271\\
81	0.29836691217747\\
82	0.421199191108479\\
83	0.662172206554038\\
84	0.240462161150455\\
85	0.315713605727069\\
86	0.638535534481614\\
87	0.718999380624762\\
88	0.777873612863486\\
89	0.272832980458975\\
90	0.364551029925225\\
91	0.394103163463619\\
92	0.416628255246637\\
93	0.400192936110889\\
94	0.831680268702745\\
95	0.903616184072511\\
96	0.699849286061203\\
97	0.466232726674374\\
98	0.726458291509927\\
99	1.02102328237711\\
100	1.22852528112014\\
};
\node[anchor=north east,font=\scriptsize,fill=white,fill opacity=0.75,text opacity=1,inner sep=1pt] at (axis cs:\xMaxAll,\yMaxNorm) {Minimax approach};
\end{axis}

\begin{axis}[
  vaxisbase,
  at={(0,0)},
  axis x line=none,
  axis y line*=right,
  ymin=\yMinB, ymax=\yMaxB,
  ylabel={$B,\hat B_t$},
  ylabel style={color=black},
]
\addplot [color=black, dashed, line width=1pt, forget plot] table[row sep=crcr]{%
0	1.95870836561711\\
1	1.95870836561711\\
99	1.95870836561711\\
};
\addplot [color=black, dashed, line width=1pt, forget plot] table[row sep=crcr]{%
0	0\\
1	0\\
99	0\\
};
\addplot [color=mycolor1, line width=1.5pt, dashdotted, forget plot] table[row sep=crcr]{%
0	0.492985075246353\\
1	1.46962223709375\\
2	1.46990513890134\\
3	1.46636739313385\\
4	1.4635545847913\\
5	1.46391055826188\\
6	1.46389836710439\\
7	1.46465541987914\\
8	1.46491980263028\\
9	1.46516967350574\\
10	1.46569929853107\\
11	1.46515116070915\\
12	1.46517436573205\\
13	1.4660728404322\\
14	1.46624653533637\\
15	1.4669782984648\\
16	1.46800692712998\\
17	1.46889660133929\\
18	1.46781570785459\\
19	1.46797432795109\\
20	1.46707760964489\\
21	1.46667302823847\\
22	1.46624992882649\\
23	1.46610058046448\\
24	1.46567120556081\\
25	1.46704640698231\\
26	1.46630939002272\\
27	1.46665300513883\\
28	1.46668983459953\\
29	1.46657339748478\\
30	1.46654906073846\\
31	1.46650815581883\\
32	1.46640243529002\\
33	1.46658725900955\\
34	1.46634245373778\\
35	1.46626701687717\\
36	1.46697766703415\\
37	1.46697547252197\\
38	1.46795059774312\\
39	1.46815259618165\\
40	1.46745085932633\\
41	1.46747028007\\
42	1.4688221771284\\
43	1.46880355747458\\
44	1.46820064354239\\
45	1.46839212378496\\
46	1.46902733240062\\
47	1.4706125083065\\
48	1.47274758008104\\
49	1.47586975151595\\
50	1.47636301443851\\
51	1.47706423895951\\
52	1.47704747934509\\
53	1.47711916848154\\
54	1.47713155814345\\
55	1.4773934230445\\
56	1.47776176662958\\
57	1.47739077209932\\
58	1.4773295865239\\
59	1.47759108128988\\
60	1.47738019213071\\
61	1.47735909237789\\
62	1.47725382957221\\
63	1.47720817112647\\
64	1.47725837357494\\
65	1.47805654147457\\
66	1.47914844117842\\
67	1.48197816788295\\
68	1.48076113903491\\
69	1.47963116888248\\
70	1.48073626913857\\
71	1.47872118837767\\
72	1.47858648343449\\
73	1.48208073304528\\
74	1.48454487063716\\
75	1.48191037148706\\
76	1.48247245325009\\
77	1.48242832262496\\
78	1.48208282397292\\
79	1.48193291514\\
80	1.4819233681232\\
81	1.48217101598868\\
82	1.48398650498607\\
83	1.48271546013\\
84	1.48352637094727\\
85	1.48344237251856\\
86	1.48452498611866\\
87	1.49052847843807\\
88	1.49084262238895\\
89	1.48887095949963\\
90	1.48881552005568\\
91	1.4908317436632\\
92	1.4907841510393\\
93	1.48852313414341\\
94	1.49159524307484\\
95	1.49713671638608\\
96	1.50041456314808\\
97	1.50401892898834\\
98	1.50767286450915\\
99	1.51808057435589\\
};
\addplot [color=mycolor1, line width=1.5pt, dashdotted, forget plot] table[row sep=crcr]{%
0	-2.95921612703286\\
1	0.833133640951125\\
2	0.833292986023249\\
3	0.831300351871045\\
4	0.829716070592278\\
5	0.829916559240458\\
6	0.829909692302239\\
7	0.830336070931263\\
8	0.830485014692997\\
9	0.830625743953764\\
10	0.830924051293604\\
11	0.830615318600726\\
12	0.830628386834341\\
13	0.831134446291695\\
14	0.831232289069088\\
15	0.831644441939492\\
16	0.832223815258241\\
17	0.832724926477833\\
18	0.832116110953602\\
19	0.832205453736389\\
20	0.831700378572327\\
21	0.831472499532992\\
22	0.831234190679451\\
23	0.831150071082853\\
24	0.83090823808349\\
25	0.83168280386994\\
26	0.831267681734229\\
27	0.831461234108719\\
28	0.831481965729876\\
29	0.831416383003974\\
30	0.831402688867983\\
31	0.831379635922226\\
32	0.831320089315701\\
33	0.83142419736477\\
34	0.831286304923266\\
35	0.831243815556825\\
36	0.831644086338918\\
37	0.831642850214104\\
38	0.832192100971245\\
39	0.832305878812751\\
40	0.831910610549727\\
41	0.83192154935918\\
42	0.832683004855131\\
43	0.832672519598675\\
44	0.83233294294538\\
45	0.832440776946484\\
46	0.832798558894991\\
47	0.833691415371223\\
48	0.834894012498483\\
49	0.836652641062242\\
50	0.836930461665125\\
51	0.837325440881359\\
52	0.837316001619345\\
53	0.83735638111653\\
54	0.83736336039031\\
55	0.837510882672558\\
56	0.837718339206635\\
57	0.837509378589691\\
58	0.837474903713817\\
59	0.837622196391263\\
60	0.837503408558074\\
61	0.837491525286351\\
62	0.837432254521932\\
63	0.837406538882995\\
64	0.837434831045312\\
65	0.837884402407489\\
66	0.838499444967928\\
67	0.840093347025849\\
68	0.83940784374059\\
69	0.838771328203104\\
70	0.839393836650089\\
71	0.838258757625678\\
72	0.838182881716344\\
73	0.840151120105067\\
74	0.841539139134056\\
75	0.840055158486582\\
76	0.840371769997714\\
77	0.840346917964422\\
78	0.840152329435238\\
79	0.840067856973223\\
80	0.840062479304642\\
81	0.840202010615677\\
82	0.841224615132505\\
83	0.840508652136652\\
84	0.840965427262391\\
85	0.840918112035909\\
86	0.841527936383773\\
87	0.84490968367957\\
88	0.845086642329374\\
89	0.843975999991609\\
90	0.843944776377207\\
91	0.845080514637354\\
92	0.845053705173835\\
93	0.843780114090335\\
94	0.84551059781053\\
95	0.848632181785316\\
96	0.850478715296417\\
97	0.852509121567865\\
98	0.854567575814909\\
99	0.860430771284317\\
};
\end{axis}

\end{tikzpicture}%

%% file: references.bib
@inproceedings{abbasi2011regret,
  author    = {Y.~{Abbasi-Yadkori} and C.~{Szepesv{\'a}ri}},
  title     = {Regret bounds for the adaptive control of linear quadratic systems},
  booktitle = {Proceedings of the 24th Annual Conference on Learning Theory},
  editor    = {S.~M.~{Kakade} and U.~{von Luxburg}},
  publisher = {PMLR},
  pages     = {1--26},
  year      = {2011},
  note       = {\bibnote{First general learning algorithm for linear-quadratic control with a regret bound growing as $\sqrt{T}$.}}
}

@article{agrawal1995sample,
  author  = {R.~{Agrawal}},
  title   = {Sample mean based index policies by {\it O}(log n) regret for the multi-armed bandit problem},
  journal = {Adv. Appl. Prob.},
  volume  = {27},
  number  = {4},
  pages   = {1054--1078},
  year    = {1995}
}

@article{anderson1977exponential,
  author  = {B.~{Anderson}},
  title   = {Exponential stability of linear equations arising in adaptive identification},
  journal = {IEEE Trans. Autom. Control},
  volume  = {22},
  number  = {1},
  pages   = {83--88},
  year    = {1977}
}

@article{annaswamy2021historical,
    author  =   {A.~M.~{Annaswamy} and A.~L.~{Fradkov}},
    year    =   {2021},
    title   =   {A historical perspective of adaptive control and learning},
    journal =   {Annu. Rev. Control},
    volume  =   {52},
    pages   =   {18--41}
}

@article{Arcari2020,
  author  = {E.~{Arcari} and L.~{Hewing} and M.~N.~{Zeilinger}},
  title   = {An approximate dynamic programming approach for dual stochastic model predictive control},
  journal = {IFAC-PapersOnLine},
  volume  = {53},
  number  = {2},
  pages   = {8105--8111},
  year    = {2020}
}

@book{astolfi2007nonlinear,
    author  =   {A.~{Astolfi} and D.~{Karagiannis} and R.~{Ortega}},
    year    =   {2007},
    title   =   {Nonlinear and Adaptive Control with Applications},
    publisher   =   {Springer}
}

@inproceedings{auer1995gambling,
  author    = {P.~{Auer} and N.~{Cesa-Bianchi} and Y.~{Freund} and R.~E.~{Schapire}},
  title     = {Gambling in a rigged casino: The adversarial multi-armed bandit problem},
  booktitle = {Proceedings of the 36th Annual Symposium on Foundations of Computer Science},
  publisher  = {IEEE Computer Society Press},
  pages      = {322--331},
  year       = {1995}
}

@article{auer2002finite,
  author  = {P.~{Auer} and N.~{Cesa-Bianchi} and P.~{Fischer}},
  title   = {Finite-time analysis of the multiarmed bandit problem},
  journal = {Mach. Learn.},
  volume  = {47},
  number  = {2},
  pages   = {235--256},
  year    = {2002}
}

@article{bai1985persistency,
  author  = {E.~W.~{Bai} and S.~S.~{Sastry}},
  title   = {Persistency of excitation, sufficient richness and parameter convergence in discrete time adaptive control},
  journal = {Syst. Control Lett.},
  volume  = {6},
  number  = {3},
  pages   = {153--163},
  year    = {1985}
}

@misc{Baltussen2026,
  author = {T.~{Baltussen} and N.~P.~{Lawrence} and A.~{Katriniok} and A.~{Mesbah} and M.~{Heemels}},
  title  = {The separation principle and the dual-certainty equivalence gap in model predictive control},
  year   = {2026},
  note   = {\url{https://arxiv.org/abs/2604.06045}}
}

@misc{Baltussen2025b,
  author = {T.~{Baltussen} and M.~{Heemels} and A.~{Katriniok}},
  title  = {Dual {MPC} for active learning of nonparametric uncertainties},
  year   = {2025},
  note   = {\url{https://arxiv.org/abs/2511.08542}}
}

@article{bar1974dual,
  author  = {Y.~{Bar-Shalom} and E.~{Tse}},
  title   = {Dual effect, certainty equivalence, and separation in stochastic control},
  journal = {IEEE Trans. Autom. Control},
  volume  = {19},
  number  = {5},
  pages   = {494--500},
  year    = {1974},
  note    = {\bibnote{Investigated important concepts for certainty-equivalence control and dual control.}}
}

@book{Basar/B95,
  author    = {T.~{Basar} and P.~{Bernhard}},
  title     = {${H}_{\infty}$-Optimal Control and Related Minimax Design Problems---{A} Dynamic Game Approach},
  publisher = {Birkhauser},
  year      = {1995}
}

@book{basar1999dynamic,
  author    = {T.~{Basar} and G.~J.~{Olsder}},
  title     = {Dynamic Noncooperative Game Theory},
  publisher = {SIAM},
  year      = {1999}
}

@book{Bellman57,
    author  =   {R.~E.~{Bellman}},
    year    =   {1957},
    title   =   {Dynamic Programming},
    publisher   =   {Princeton Univ. Press}
}

@inproceedings{bencherki2025,
  author    = {F.~{Bencherki} and A.~{Rantzer}},
  title     =   {Adaptive control of positive systems with application to learning {SSP}},
  booktitle = {Proceedings of the 7th Annual Learning for Dynamics and Control Conference},
  editor     = {N.~{Ozay} and L.~{Balzano} and D.~{Panagou} and A.~{Abate}},
  pages      = {660--672},
  publisher  = {PMLR},
  year       = {2025}
}

@inproceedings{Berberich2023,
  author    = {J.~{Berberich} and A.~{Ianelli} and A.~{Padoan} and J.~{Coulson} and F.~{{D\"{o}rfler}} and F.~{{Allg\"{o}wer}}},
  title     = {A quantitative and constructive proof of {Willems}' fundamental lemma and its implications},
  booktitle = {American Control Conference},
  publisher = {IEEE},
  pages     = {4155--4160},
  year      = {2023}
}

@article{bernhardsson1989dual,
    author  =   {B.~{Bernhardsson}},
    year    =   {1989},
    title   =   {Dual control of a first-order system with two possible gains},
    journal =   {Int. J. Adapt. Control Signal Process.},
    volume  =   {3},
    pages   =   {15--22}
}

@book{bertsekas2007dynamic,
  author    = {D.~P.~{Bertsekas} and et~al.},
  title     = {Dynamic Programming and Optimal Control, volume {II}},
  publisher = {Athena Scientific},
  edition   = {3rd},
  year      = {2007}
}

@book{Bertsekas2012,
  author    = {D.~P.~{Berstekas}},
  title     = {Approximate Dynamic Programming},
  publisher = {Athena Scientific},
  edition   = {4th},
  year      = {2012}
}

@article{campi1998adaptive,
  author  = {M.~C.~{Campi} and P.~{Kumar}},
  title   = {Adaptive linear quadratic {Gaussian} control: {The} cost-biased approach revisited},
  journal = {SIAM J. Control Optim.},
  volume  = {36},
  number  = {6},
  pages   = {1890--1907},
  year    = {1998}
}

@inproceedings{cassel2020logarithmic,
  author    = {A.~{Cassel} and A.~{Cohen} and T.~{Koren}},
  title     = {Logarithmic regret for learning linear quadratic regulators efficiently},
  booktitle = {Proceedings of the 37th International Conference on Machine Learning},
  editor    = {H.~{Daum{\'e} III} and A.~{Singh}},
  pages     = {1328--1337},
  publisher = {PMLR},
  year      = {2020}
}

@article{charnov1976optimal,
    author  =   {E.~L.~{Charnov}},
    title   =   {Optimal foraging, the marginal value theorem},
    journal =   {Theor. Popul. Biol.},
    volume  =   {9},
    number  =   {2},
    pages   =   {129--136},
    year    =   {1976}
}

@book{chen2012identification,
    author  =   {H.~F.~{Chen} and L.~{Guo}},
    year    =   {2012},
    title   =   {Identification and Stochastic Adaptive Control},
    publisher   =   {Springer}
}

@article{cohen2007should,
    author  =   {J.~D.~{Cohen} and S.~M.~{McClure}},
    title   =   {Should {I} stay or should {I} go? {How} the human brain manages the trade-off between exploitation and exploration},
    journal =   {Philos. Trans. R. Soc. B Biol. Sci.},
    volume  =   {362},
    number  =   {1481},
    pages   =   {933--942},
    year    =   {2007}
}

@inproceedings{Cusumano/P88,
  author    = {S.~{Cusumano} and K.~{Poolla}},
  title     = {Nonlinear feedback vs. linear feedback for robust stabilization},
  booktitle = {IEEE Conference on Decision and Control},
  publisher = {IEEE},
  pages     = {1776--1780},
  year      = {1988}
}

@inproceedings{Didinsky1994,
  author    = {G.~{Didinsky} and T.~{{Ba\c{s}ar}}},
  title     = {Minimax adaptive control of uncertain plants},
  booktitle = {IEEE Conference on Decision and Control},
  publisher = {IEEE},
  pages     = {2839--2844},
  year      = {1994}
}

@article{fan1953minimax,
  author  = {K.~{Fan}},
  title   = {Minimax theorems},
  journal = {Proc. Natl. Acad. Sci. USA},
  volume  = {39},
  number  = {1},
  pages   = {42--47},
  year    = {1953}
}

@article{feldbaum1960dual,
    author  =   {A.~A.~{Feldbaum}},
    title   =   {Dual control theory {I}},
    journal =   {Avtomat. i Telemekh},
    volume  =   {21},
    number  =   {9},
    pages   =   {1240--1249},
    year    =   {1960},
    note    =   {\bibnote{Pioneered the concept of dual control and recognized the need to balance learning and regulation.}}
}

@article{feldbaum1963dual,
    author  =   {A.~A.~{Feldbaum}},
    title   =   {Dual control theory problems},
    journal =   {IFAC Proc. Vol.},
    volume  =   {1},
    number  =   {2},
    pages   =   {541--550},
    year    =   {1963}
}

@article{gittins1974dynamic,
    author  =   {J.~{Gittins} and D.~{Jones}},
    year    =   {1974},
    title   =   {A dynamic allocation index for the sequential design of experiments},
    journal =   {Prog. Stat.},
    volume  =   {1},
    number  =   {1},
    pages   =   {241--266},
    note    =   {\bibnote{Solved the multi-armed bandit problem in a Bayesian setting.}}
}

@book{goodwin2014adaptive,
    author  =   {G.~C.~{Goodwin} and K.~S.~{Sin}},
    year    =   {2014},
    title   =   {Adaptive Filtering Prediction and Control},
    publisher   =   {Courier}
}

@article{guo1995convergence,
  author  = {L.~{Guo}},
  title   = {Convergence and logarithm laws of self-tuning regulators},
  journal = {Automatica},
  volume  = {31},
  number  = {3},
  pages   = {435--450},
  year    = {1995}
}

@article{guo1996self,
  author  = {L.~{Guo}},
  title   = {Self-convergence of weighted least-squares with applications to stochastic adaptive control},
  journal = {IEEE Trans. Automatic Control},
  volume  = {41},
  number  = {1},
  pages   = {79--89},
  year    = {1996}
}

@article{gurpegui2023minimax,
  author  = {A.~{Gurpegui} and E.~{Tegling} and A.~{Rantzer}},
  title   = {Minimax linear optimal control of positive systems},
  journal   =   {IEEE Control Syst. Lett.},
  volume  = {7},
  pages   = {3920--3925},
  year    = {2023}
}

@book{holland1992adaptation,
    author  =   {J.~H.~{Holland}},
    title   =   {Adaptation in Natural and Artificial Systems: {An} Introductory Analysis with Applications to Biology, Control, and Artificial Intelligence},
    year    =   {1992},
    publisher   =   {MIT press}
}

@article{Ioannou1988,
  author  = {P.~{Ioannou} and J.~{Sun}},
  title   = {Theory and design of robust direct and indirect adaptive-control schemes},
  journal = {Int. J. Control},
  volume  = {1988},
  number  = {3},
  pages   = {775--813},
  year    = {1988}
}

@inproceedings{jedra2022minimal,
  author    = {Y.~{Jedra} and A.~{Proutiere}},
  title     = {Minimal expected regret in linear quadratic control},
  booktitle = {Proceedings of the 25th International Conference on Artificial Intelligence and Statistics},
  editor    = {G.~{Camps-Valls} and F.~J.~R.~{Ruiz} and I.~{Valera}},
  pages     = {10234--10321},
  publisher = {PMLR},
  year      = {2022}
}

@article{jung2021systematic,
    author  =   {D.~{Jung} and Y.~{Choi}},
    title   =   {Systematic review of machine learning applications in mining: {Exploration}, xploitation, and reclamation.},
    journal =   {Minerals},
    volume  =   {11},
    number  =   {2},
    pages   =   {148},
    year    =   {2021}
}

@article{kelly1981multi,
  author  = {F.~P.~{Kelly}},
  title   = {Multi-armed bandits with discount factor near one: {The} {Bernoulli} case},
  journal = {Ann. Stat.},
  volume  = {9},
  number  = {5},
  pages   = {987--1001},
  year    = {1981}
}

@inproceedings{kjellqvist2022learning,
  author    = {O.~{Kjellqvist} and A.~{Rantzer}},
  title     = {Learning-enabled robust control with noisy measurements},
  booktitle = {Proceedings of the 4th Annual Learning for Dynamics and Control Conference},
  editor    = {R.~{Firoozi} and N.~{Mehr} and E.~{Yel} and R.~{Antonova} and J.~{Bohg} and M.~{Schwager} and M.~{Kochenderfer}},
  pages     = {86--96},
  publisher = {PMLR},
  year      = {2022}
}

@inproceedings{kjellqvist2022minimax,
  author    = {O.~{Kjellqvist} and A.~{Rantzer}},
  title     = {Minimax adaptive estimation for finite sets of linear systems},
  booktitle = {American Control Conference},
  publisher = {IEEE},
  pages     = {260--265},
  year      = {2022}
}

@article{kulcsar1996dual,
    author  =   {C.~{Kulcs{\'a}r} and L.~{Pronzato} and E.~{Walter}},
    year    =   {1996},
    title   =   {Dual control of linearly parameterised models via prediction of posterior densities},
    journal =   {Eur. J. Control},
    volume  =   {2},
    number  =   {2},
    pages   =   {135--143}
}

@article{kumar1985survey,
    author  =   {P.~R.~{Kumar}},
    year    =   {1985},
    title   =   {A survey of some results in stochastic adaptive control},
    journal =   {SIAM J. Control Optim.},
    volume  =   {23},
    number  =   {3},
    pages   =   {329--380},
}

@article{lai1978adaptive,
  author  = {T.~{Lai} and H.~{Robbins}},
  title   = {Adaptive design in regression and control},
  journal = {Proc. Natl. Acad. Sci. USA},
  volume  = {75},
  number  = {2},
  pages   = {586--587},
  year    = {1978}
}

@article{lai1985asymptotically,
  author  = {T.~L.~{Lai} and H.~{Robbins}},
  title   = {Asymptotically efficient adaptive allocation rules},
  journal = {Adv. Appl. Math.},
  volume  = {6},
  number  = {1},
  pages   = {4--22},
  year    = {1985},
  note    = {\bibnote{Established the logarithmic lower bound on achievable regret for bandits from a non-Bayesian perspective.}}
}

@article{lai1987asymptotically,
  author  = {T.~L.~{Lai} and C.~Z.~{Wei}},
  title   = {Asymptotically efficient self-tuning regulators},
  journal = {SIAM J. Control Optim.},
  volume  = {25},
  number  = {2},
  pages   = {466--481},
  year    = {1987},
  note    = {\bibnote{Pioneered logarithmic regret growth in the context of self-tuning control.}}
}

@article{lavie2010exploration,
    author  =   {D.~{Lavie} and U.~{Stettner} and M.~L.~{Tushman}},
    title   =   {Exploration and exploitation within and across organizations},
    journal =   {Acad. Manag. Ann.},
    volume  =   {4},
    number  =   {1},
    pages   =   {109--155},
    year    =   {2010}
}

@inproceedings{lee2024nonasymptotic,
  author    = {B.~{Lee} and A.~{Rantzer} and N.~{Matni}},
  title     = {Nonasymptotic regret analysis of adaptive linear quadratic control with model misspecification},
  booktitle = {Proceedings of the 6th Annual Learning for Dynamics and Control Conference},
  editor    = {A.~{Abate} and M.~{Cannon} and K.~{Margellos} and A.~{Papachristodoulou}},
  pages     = {980--992},
  publisher = {PMLR},
  year      = {2024}
}

@article{ljung1975counterexamples,
  author  = {L.~{Ljung} and T.~{{S{\"o}derstr{\"o}m}} and I.~{Gustavsson}},
  title   = {Counterexamples to general convergence of a commonly used recursive identification method},
  journal = {IEEE Trans. Autom. Control},
  volume  = {20},
  number  = {5},
  pages   = {643--652},
  year    = {1975}
}

@article{ljung77on,
  author  = {L.~{Ljung}},
  title   = {On positive real transfer function and convergence of some recursive schemes},
  journal = {IEEE Trans. Autom. Control},
  volume  = {AC-22},
  pages   = {539--551},
  year    = {1977}
}

@article{macarthur1966optimal,
    author  =   {R.~H.~{MacArthur} and E.~R.~{Pianka}},
    title   =   {On optimal use of a patchy environment},
    journal =   {Am. Nat.},
    volume  =   {100},
    number  =   {916},
    pages   =   {603--609},
    year   =   {1966}
}

@article{Markovsky2008,
  author  = {I.~{Markovsky} and P.~{Rapisarda}},
  title   = {Data-driven simulation and control},
  journal = {Int. J. Control},
  volume  = {81},
  number  = {12},
  pages   = {1946--1959},
  year    = {2008}
}

@misc{Meijer2025,
  author = {T.~J.~{Meijer} and K.~J.~A.~{Scheres} and S.~A.~N.~{Nouwens} and V.~S.~{Dolk} and W.~P.~M.~H.~{Heemels}},
  title  = {From a frequency-domain {Willems'} lemma to data-driven predictive control},
  year   = {2025},
  note   = {\url{https://arxiv.org/abs/2501.19390}}
}

@article{mesbah2018stochastic,
  author  = {A.~{Mesbah}},
  title   = {Stochastic model predictive control with active uncertainty learning: {A} survey on dual control},
  journal = {Annu. Rev. Control},
  volume  = {45},
  pages   = {107--117},
  year    = {2018}
}

@techreport{Megretski/R03,
  author      = {A.~{Megretski} and A.~{Rantzer}},
  title       = {Bounds on the optimal $l_2$ gain in adaptive control of a first order linear system},
  institution = {Institut Mittag-Leffler},
  number      = {41},
  year        = {2003}
}

@inproceedings{2004megretskinonlinear,
  author    = {A.~{Megretski}},
  title     = {A nonlinear dynamical game interpretation of adaptive $\ell_2$ control: {Performance} limitations and suboptimal controllers},
  booktitle = {Proceedings of the 16th International Symposium on Mathematical Theory of Networks and Systems},
  year      = {2004}
}

@article{Messerer2022,
  author  = {F.~{Messerer} and K.~{{Baumg\"{a}rtner}} and M.~{Diehl}},
  title   = {A dual-control effect preserving formulation for nonlinear output-feedback stochastic model predictive control with constraints},
  journal = {IEEE Control Syst. Lett.},
  volume  = {7},
  pages   = {1171--1176},
  year    = {2023}
}

@article{morgan1977stability,
  author  = {A.~{Morgan} and K.~{Narendra}},
  title   = {On the stability of nonautonomous differential equations {$\dot{x}=\left[A+B(t)\right]x$}, with skew symmetric matrix {$B(t)$}},
  journal = {SIAM J. Control Optim.},
  volume  = {15},
  number  = {1},
  pages   = {163--176},
  year    = {1977}
}

@misc{Mulagaleti2026,
    author  =   {S.~K.~{Mulagaleti} and A.~{Bemporad}},
    year    =   {2026},
    title   =   {Dual {MPC} for quasi-linear parameter varying systems},
    note    =   {Preprint: \url{https://arxiv.org/abs/2603.29445}}
}

@article{neumann1928theorie,
  author  = {J.~v.~{Neumann}},
  title   = {Zur theorie der gesellschaftsspiele},
  journal = {Mathematische annalen},
  volume  = {100},
  number  = {1},
  pages   = {295--320},
  year    = {1928}
}

@inproceedings{orlov2018adaptive,
  author    = {Y.~{Orlov} and A.~{Rantzer} and L.~T.~{Aguilar}},
  title     = {Adaptive ${H}_\infty$ synthesis for linear systems with uncertain parameters},
  booktitle = {IEEE Conference on Decision and Control},
  publisher = {IEEE},
  pages     = {5512--5517},
  year      = {2018}
}

@inproceedings{ouyang2017control,
  author    = {Y.~{Ouyang} and M.~{Gagrani} and R.~{Jain}},
  title     = {Control of unknown linear systems with {Thompson} sampling},
  booktitle = {Annual Allerton Conference on Communication, Control, and Computing},
  pages     = {1198--1205},
  publisher = {IEEE},
  year      = {2017}
}

@article{Pan1998,
  author  = {Z.~{Pan} and T.~{{Ba\c{s}ar}}},
  title   = {Adaptive controller design for tracking and disturbance attenuation in parametric strict-feedback nonlinear systems},
  journal = {IEEE Trans. Autom. Control},
  volume  = {43},
  number  = {8},
  pages   = {1066--1083},
  year    = {1998}
}

@article{Parsi2023a,
  author  = {A.~{Parsi} and A.~{Iannelli} and R.~S.~{Smith}},
  title   = {An explicit dual control approach for constrained reference tracking of uncertain linear systems},
  journal = {IEEE Trans. Autom. Control},
  volume  = {68},
  number  = {5},
  pages   = {2652--2666},
  year    = {2023}
}

@article{Parsi2023b,
  author  = {A.~{Parsi} and D.~{Liu} and A.~{Iannelli} and R.~S.~{Smith}},
  title   = {Dual adaptive {MPC} using an exact set-membership reformulation},
  journal = {IFAC-PapersOnLine},
  volume  = {56},
  number  = {2},
  pages   = {8457--8463},
  year    = {2023}
}

@article{rantzer2020ifac,
  author  = {A.~{Rantzer}},
  title   = {Minimax adaptive control for state matrix with unknown sign},
  journal = {IFAC-PapersOnLine},
  volume  = {53},
  number  = {2},
  pages   = {58--62},
  year    = {2020},
  note    = {Corrected version avaliable from arXiv:1912.03550}
}

@inproceedings{rantzer2021minimax,
  author    = {A.~{Rantzer}},
  title     = {Minimax adaptive control for a finite set of linear systems},
  booktitle = {Proceedings of the 3rd Conference on Learning for Dynamics and Control},
  editor    = {A.~{Jadbabaie} and J.~{Lygeros} and G.~J.~{Pappas} and P.~A.~{Parrilo} and B.~{Recht} and C.~J.~{Tomlin} and M.~N.~{Zeilinger}},
  pages     = {893--904},
  publisher = {PMLR},
  year      = {2021},
  note      = {\bibnote{The minimax dual control problem is reformulated as a standard zero-sum dynamic game.}}
}

@article{RantzerValcher2021,
  author  = {A.~{Rantzer} and M.~E.~{Valcher}},
  title   = {Scalable control of positive systems},
  journal = {Annu. Rev. Control Robot. Auton. Syst.},
  volume  = {4},
  pages   = {319--341},
  year    = {2020}
}

@inproceedings{rantzer2022explicit,
  author    = {A.~{Rantzer}},
  title     = {Explicit solution to {B}ellman equation for positive systems with linear cost},
  booktitle = {IEEE Conference on Decision and Control},
  publisher = {IEEE},
  pages     = {6154--6155},
  year      = {2022}
}

@inproceedings{Rantzer2025acc,
  author    = {A.~{Rantzer}},
  title     = {On minimax optimal dual control for fully actuated systems},
  booktitle = {American Control Conference},
  publisher = {IEEE},
  pages     = {3993--3997},
  year      = {2025}
}

@article{rantzer2026minimax,
  author  = {A.~{Rantzer}},
  title   = {Minimax optimal dual control--{The} single input case},
  journal = {arXiv preprint arXiv:2604.18550},
  year    = {2026}
}

@article{recht2019tour,
  author  = {B.~{Recht}},
  title   = {A tour of reinforcement learning: {The} view from continuous control},
  journal = {Annu. Rev. Control Robot. Auton. Syst.},
  volume  = {2},
  number  = {1},
  pages   = {253--279},
  year    = {2019}
}

@book{sastry2011adaptive,
    author  =   {S.~{Sastry} and M.~{Bodson}},
    year    =   {2011},
    title   =   {Adaptive Control: {Stability}, Convergence and Robustness},
    publisher   =   {Courier}
}

@inproceedings{simchowitz2020naive,
  author    = {M.~{Simchowitz} and D.~{Foster}},
  title     = {Naive exploration is optimal for online {LQR}},
  booktitle = {Proceedings of the 37th International Conference on Machine Learning},
  editor    = {H.~{Daum{\'e} III} and A.~{Singh}},
  pages     = {8937--8948},
  publisher = {PMLR},
  year      = {2020}
}

@article{soderstrom1975identifiability,
  author  = {T.~{{S{\"o}derstr{\"o}m}} and I.~{Gustavsson} and L.~{Ljung}},
  title   = {Identifiability conditions for linear systems operating in closed loop},
  journal = {Int. J. Control},
  volume  = {21},
  number  = {2},
  pages   = {243--255},
  year    = {1975}
}

@article{sternby1976simple,
    author  =   {J.~{Sternby}},
    year    =   {1976},
    title   =   {A simple dual control problem with an analytical solution},
    journal =   {IEEE Trans. Autom. Control},
    volume  =   {21},
    number  =   {6},
    pages   =   {840--844}
}

@book{sutton1998reinforcement,
    author  =   {R.~S.~{Sutton} and A.~G.~{Barto}},
    year    =   {1998},
    title   =   {Reinforcement learning: {An} introduction},
    publisher   =   {MIT press}
}

@article{thompson1933likelihood,
    author  =   {W.~R.~{Thompson}},
    title   =   {On the likelihood that one unknown probability exceeds another in view of the evidence of two samples},
    journal =   {Biometrika},
    volume  =   {25},
    number  =   {3/4},
    pages   =   {285--294},    
    year    =   {1933}
}

@article{tsiamis2023statistical,
    author  =   {A.~{Tsiamis} and I.~{Ziemann} and N.~{Matni} and G.~J.~{Pappas}},
    year    =   {2023},
    title   =   {Statistical learning theory for control: {A} finite-sample perspective},
    journal =   {IEEE Control Syst. Mag.},
    volume  =   {43},
    number  =   {6},
    pages   =   {67--97}
}

@book{vanOverschee1996,
  author    = {P.~{{van Overschee}} and B.~{{De Moor}}},
  title     = {Subspace identification for linear systems: {Theory}, implementation, applications},
  publisher = {Kluwer Academic Publishers},
  year      = {1996}
}

@article{vanWaarde2020c,
  author  = {H.~J.~{{van Waarde}} and J.~{Eising} and H.~L.~{Trentelmann} and M.~K.~{Camlibel}},
  title   = {Data informativity: {A} new perspective on data-driven analysis and control},
  journal = {IEEE Trans. Autom. Control},
  volume  = {65},
  number  = {11},
  pages   = {4753--4768},
  year    = {2020}
}

@inproceedings{vinnicombe2004examples,
  author    = {G.~{Vinnicombe}},
  title     = {Examples and counterexamples in finite {L2}-gain adaptive control},
  booktitle = {Proceedings of the 16th International Symposium on Mathematical Theory of Networks and Systems},
  year      = {2004},
  note       = {\bibnote{Established fundamental performance limitations and benchmarks for minimax dual control.}}
}

@article{Wang2015,
  author  = {Y.~{Wang} and B.~{O'Donoghue} and S.~{Boyd}},
  title   = {Approximate dynamic programming via iterated {Bellman} inequalities},
  journal = {Int. J. Robust Nonlinear Control},
  volume  = {25},
  number  = {10},
  pages   = {1472--1496},
  year    = {2015}
}

@article{whittle1979discussion,
    author  =   {P.~{Whittle}},
    title   =   {Discussion of {D}r {G}ittins’ paper},
    journal =   {J. R. Stat. Soc.},
    volume  =   {41},
    number  =   {2},
    pages   =   {164--177},
    year    =   {1979}
}

@article{wittenmark1995adaptive,
    author  =   {B.~{Wittenmark}},
    year    =   {1995},
    title   =   {Adaptive dual control methods: {An} overview},
    journal =   {IFAC Proc. Vol.},
    volume  =   {28},
    number  =   {13},
    pages   =   {67--72}
}

@article{Zames1981,
  author  = {G.~{Zames}},
  title   = {Feedback and optimal sensitivity: {Model} reference transformations, multiplicative seminorms, and approximate inverses},
  journal = {IEEE Trans. Autom. Control},
  volume  = {26},
  number  = {2},
  pages   = {301--320},
  year    = {1981}
}

@book{zhou+96,
  author    = {K.~{Zhou} and J.~{Doyle} and K.~{Glover}},
  title     = {Robust and Optimal Control},
  publisher = {Prentice-Hall},
  year      = {1996}
}

@article{ziemann2020phase,
  author  = {I.~{Ziemann} and H.~{Sandberg}},
  title   = {On a phase transition of regret in linear quadratic control: {The} memoryless case},
  journal = {IEEE Control Syst. Lett.},
  volume  = {5},
  number  = {2},
  pages   = {695--700},
  year    = {2020}
}

@article{ziemann2025regret,
  author  = {I.~{Ziemann} and H.~{Sandberg}},
  title   = {Regret lower bounds for learning linear quadratic {Gaussian} systems},
  journal = {IEEE Trans. Autom. Control},
  volume  = {70},
  number  = {1},
  pages   = {159--173},
  year    = {2025}
}

@article{aastrom1965optimal,
    author  =   {K.~J.~{{\AA}str{\"o}m}},
    year    =   {1965},
    title   =   {Optimal control of {Markov} processes with incomplete state information},
    journal =   {J. Math. Anal. Appl.},
    volume  =   {10},
    number  =   {1},
    pages   =   {174--205}
}

@article{aastrom1965numerical,
  author  = {K.~J.~{{\AA}str{\"o}m} and T.~{Bohlin}},
  title   = {Numerical identification of linear dynamic systems from normal operating records},
  journal = {IFAC Proc. Vol.},
  volume  = {2},
  number  = {2},
  pages   = {96--111},
  year    = {1965},
  note    = {\bibnote{Introduced persistence of excitation conditions for consistent parameter estimation in adaptive systems.}}
}

@article{ast+wit73,
  author  = {K.~J.~{{\AA}str{\"o}m} and B.~{Wittenmark}},
  title   = {On self-tuning regulators},
  journal = {Automatica},
  volume  = {9},
  number  = {2},
  pages   = {185--199},
  year    = {1973},
  note    = {\bibnote{Seminal paper introducing the self-tuning regulator, a cornerstone of adaptive control.}}
}

@article{86aastrom+,
    author  =   {K.~J.~{{\AA}str{\"o}m} and A.~{Helmersson}},
    year    =   {1986},
    title   =   {Dual control of an integrator with unknown gain},
    journal =   {Comput. Math. Appl.},
    volume  =   {12},
    number  =   {6A},
    pages   =   {653--662}
}

@book{aastrom2013adaptive,
    author  =   {K.~J.~{{\AA}str{\"o}m} and B.~{Wittenmark}},
    year    =   {2013},
    title   =   {Adaptive Control},
    publisher   =   {Courier}
}

@phdthesis{kjellqvist2024minimax,
  title={Minimax Adaptive Control and Estimation},
  author={Kjellqvist, Olle},
  year={2024},
  school={Lund University}
}
